\documentclass[12pt,hyperref]{article}
\usepackage{graphicx} 
\usepackage{float} 
\usepackage{amssymb}
\usepackage{amsmath}
\usepackage{mathtools}
\usepackage[thmmarks,amsmath]{ntheorem} 
\usepackage{lineno}

\usepackage{algorithm}
\usepackage{algpseudocode}
\usepackage{amsmath}
\renewcommand{\algorithmicrequire}{\textbf{Input:}}  
\renewcommand{\algorithmicensure}{\textbf{Output:}} 

\usepackage{bm} 
\usepackage{extarrows}
\usepackage{mathrsfs} 
\usepackage{graphicx}
\usepackage{subcaption}
\usepackage{textcomp}
\usepackage{algorithm,algpseudocode,float}
\usepackage{lipsum}


\usepackage{booktabs} 
{
    \theoremstyle{nonumberplain}
    \theoremheaderfont{\bfseries}
    \theorembodyfont{\normalfont}
    \theoremsymbol{\mbox{$\Box$}}
    \newtheorem{proof}{Proof}
}

\usepackage{theorem}
\newtheorem{theorem}{Theorem}[section]

\newtheorem{lemma}{Lemma}[section]

\newtheorem{fact}{Fact}[section]
\newtheorem{claim}{Claim}[section]
\newtheorem{property}{Property}[section]
\newtheorem{question}{Question}[section]
\newtheorem{problem}{Problem}[section]
\newtheorem{conjecture}{Conjecture}[section]
{
    \theoremheaderfont{\bfseries}
    \theorembodyfont{\normalfont}
    
}
\newcommand{\RNum}[1]{\uppercase\expandafter{\romannumeral #1\relax}}

\renewcommand{\figurename}{Fig.}
\graphicspath{{figure/}}                                 
\usepackage[a4paper]{geometry}                        
\begin{document}

\title{Solutions to the linkage conjecture in tournaments \thanks{The author's work is supported by National Natural Science Foundation of China (No.12071260)}}

\author{Jia Zhou, Jin Yan\thanks{Corresponding author. E-mail adress: yanj@sdu.edu.cn}  \unskip\\[2mm]
School of Mathematics, Shandong University, Jinan 250100, China}

\date{}
\maketitle


\begin{abstract}
A digraph $D$ is $k$-linked if for every $2k$-tuple $ x_1,\ldots , x_k, y_1, \ldots , y_k$ of distinct vertices in $D$, there exist $k$ pairwise vertex-disjoint paths $P_1,\ldots, P_k$ such that $P_i$ starts at $x_i$ and ends at $y_i$, $i\in [k]$. In 2015, Pokrovskiy conjectured that there exists a function $g(k)$ such that every $2k$-connected tournament with minimum in-degree and minimum out-degree at least $g(k)$ is $k$-linked in [J. Comb. Theory, Ser. B 115 (2015) 339--347]. In this paper, we disprove this conjecture by constructing a family of counterexamples. The counterexamples also provide a negative answer to the question raised by Gir\~{a}o, Popielarz, Snyder in [Combinatorica 41 (2021) 815--837]. Further, we prove that every $(2k+1)$-connected semicomplete digraph $D$ with minimum out-degree at least $ 10^7k^{4}$ is $k$-linked, which refines and generalizes the early result of Gir\~{a}o, Popielarz, Snyder. 
\end{abstract}		

\vspace{1ex}
{\noindent\small{\bf Keywords: } connectivity; linkage; tournaments; semicomplete digraphs }
\vspace{1ex}

{\noindent\small{\bf AMS subject classifications.} 05C20, 05C38, 05C40}

\section{Introduction}
	Connectivity is one of the most fundamental concepts in graphs. Let $k$ be a positive integer. A (di)graph is \emph{$k$-connected}, if it has at least $k+1$ vertices and after deleting any set $S$ of at most $k-1$ vertices, there exists a path from $x$ to $y$, for any pair of distinct vertices $x$ and $y$. Menger's Theorem, as a classical result of connectivity, gives a necessary and sufficient condition for the existence of vertex-disjoint paths from one vertex set to the other. However, on many occasions, we still need to specify the initial and terminal vertices of each path. This yields the definition of linkage. A digraph (resp. graph) $D$ is \emph{$k$-linked} if for every $2k$-tuple $ x_1,\ldots , x_k, y_1, \ldots , y_k$ of distinct vertices in $D$, there exist $k$ pairwise vertex-disjoint paths $P_1,\ldots, P_k$ such that $P_i$ starts at $x_i$ and ends at $y_i$, $i\in [k]$. Clearly, every $k$-linked graph is $k$-connected. The converse, however, seems far from true. The $k$-linkage is clearly a much stronger property than $k$-connectivity. Is there a function $f(k)$ such that a $f(k)$-connected graph is $k$-linked? Ballob\'{a}s and Thomason \cite{Bollobas(1996)} first gave a linear function $f(k)$ with $f(k)\leq 22k$. Later, Thomas and Wollan \cite{Thomas(2005)} showed that every $2k$-connected graph with minimum degree $10k$ is $k$-linked, which is the best bound currently.

Naturally, Thomassen \cite{Thomassen(1980)} conjectured that there exists a function $f(k)$ such that every $f(k)$-connected  digraph is $k$-linked. Indeed, in 1991, Thomassen \cite{Thomassen(1991)} also disproved this conjecture. But the situation in tournaments is completely different. \emph{Tournaments} are digraphs obtained by orienting each edge of a complete graph in exactly one direction. K\"{u}hn, Lapinskas, Osthus, and Patel in \cite{Khn(2014)} concluded that any $10^4k \text{log} k$-connected tournament is $k$-linked, and they conjectured that the bound of connectivity is linear in $k$. Subsequently, in 2015, this conjecture was confirmed by Pokrovskiy \cite{Pokrovskiy(2015)} by showing that every $452k$-connected tournament is $k$-linked. Inspired by the result in \cite{Thomas(2005)} on graphs, he conjectured that a connectivity of $2k$ supplemented with certain minimum semidegree would be sufficient and thus proposed the following conjecture. Let $\delta^+(D)$ (resp., $\delta^-(D)$) denote the \emph{minimum out-degree} (resp., the \emph{minimum in-degree}) of $D$, and let $\delta^0(D)=\min \{\delta^+(D), \delta^-(D)\} $ denote the \emph{ minimum semidegree}.

\begin{conjecture}\label{conj1}
	\cite{Pokrovskiy(2015)} For every $k\in \mathbb{N}$, there exists an integer $g(k)$ such that every $2k$-connected tournament with $\delta^0(D)\geq g(k)$ is $k$-linked.
\end{conjecture}	

We disprove Conjecture \ref{conj1} by constructing counterexamples with arbitrarily large minimum semidegree.
	\begin{theorem}\label{the1}
		There exist a family of $2k$-connected tournaments of order $n$ with minimum semidegree at least $\lfloor \dfrac{n-2k}{2k+2} \rfloor $, which are not $k$-linked.
	\end{theorem}

 Gir\~{a}o, Popielarz, and Snyder \cite{Girao(2021)} showed that  the connectivity of $2k + 1$ is sufficient in tournaments with large minimum out-degree.
\begin{theorem}\label{the2}
\cite{Girao(2021)} Every $(2k + 1)$-connected tournament with minimum out-degree at least $Ck^{31}$ is $k$-linked for some constant $C$.
\end{theorem}

They also raised the following question.
\begin{question}\label{que1}
Does Theorem \ref{the2} still hold if we replace $2k+1$ by $2k$?
\end{question}

Theorem \ref{the1} indicates that a $2k$-connected tournament with the semidegree associated with \(n\) is not $k$-linked. This also means that the connectivity of $2k+1$ is tight if the lower bound of the semidegree is a function of $k$, which gives a negative answer to Question \ref{que1}. There are some related results regarding Conjecture \ref{conj1}, see \cite{bang(2021),Snyder(2019),Meng(2021)}. Gir\~{a}o, Popielarz, and Snyder \cite{Girao(2021)} also constructed an infinite family of $(2.5k- 1)$-connected tournaments that are not $k$-linked, yielding that the condition on large minimum out-degree condition is necessary when the connectivity is small. Inspired by these results, we prove the following theorem, which also provides a better lower bound of out-degree than that of \cite{Girao(2021)}. Meanwhile, it can be proved to hold on \emph{semicomplete digraphs}, the digraphs obtained by replacing every edge of a complete graph in at least one direction. Clearly, the class of semicomplete digraphs is a generalization of tournaments.


\begin{theorem}\label{cor1}
Every $(2k+1)$-connected semicomplete digraph $D$ with $\delta^+(D)\geq  10^7k^{4}$ is $k$-linked.
\end{theorem}

Bang-Jensen and Jord\'{a}n \cite{Bang(2010)} conjectured that there exists a spanning $k$-connected tournament in a $(2k - 1)$-connected semicomplete digraph, which is still open for $k\geq 4$. This conjecture indicates that we could need relatively large connectivity to ensure that a semicomplete digraph is $k$-linked, approximately twice that of a tournament.



	The rest of the paper is organized as follows. Section 2 provides notation that be utilized throughout the article. The counterexamples in Theorem \ref{the1} are given in Section 3, and in Section 4, we first show some preliminary lemmas, then we give the proof of Theorem \ref{cor1}. We use the basic strategy in \cite{Girao(2021)}. The main novelty in our proof is that we devise an algorithm (Lemma \ref{lemma4}) and also design a switching path operation, which leads to the improved bound in Theorem \ref{cor1}.

\section{Basic terminology}
Let $ \mathbb{N} $ be the set of natural numbers and for an integer $k$, we denote $[k]= \{1,\ldots,k\}$, and we define $[k,k+i]=\{k,k+1,\ldots,k+i\}$. The digraphs considered in this paper are always simple (without loops and multiple arcs). Let $D$ be a digraph with vertex set $V(D)$ and arc set $A(D)$. We write $|D|$ for the number of vertices in $D$. Given two vertices $x, y \in V (D)$, we write $xy$ for the arc from $x$ to $y$. For a vertex $x$ of a digraph $D$, we write $N^+_D (x) = \{y\ |\ xy\in A(D)\}$ (resp., $N^-_D (x) = \{y\ |\  yx\in A(D)\}$) to be the \emph{out-neighbourhood} (resp., \emph{in-neighbourhood}) and $d^+_D (x) = |N^+_D (x)|$ (resp., $d^-_D (x) = |N^-_D (x)|$) to be the \emph{out-degree} (resp., \emph{in-degree}) of $x$. Given a set $X\subseteq V(D)$, the subgraph of $D$ induced by $X$ is denoted by $D[ X ] $ and let $D - X$ or $D\setminus X$ denote the digraph obtained from $D$ by deleting $X$ and all arcs incident with $X$. We also say that $X$ is an $|X|$-element set (or just an $|X|$-set). Given a digraph $D$ and sets $A,B\subseteq V (D)$, we write $A\rightarrow B$ if every arc of $D$ between $A$ and $B$ is directed from $A$ to $B$. 

A \emph{matching} in a digraph is a set of arcs that are from $A$ to $B$ without common vertices. A matching $M$ is \emph{maximum} if there is no arc in $D\setminus M$. If $P=x_1\cdots x_t$ is a directed path of $D$, we refer to $x_1(x_t)$ as the \emph{initial}(\emph{terminal}) vertex of $P$, written as $(x_1,x_t)$-path. The path $P$ is \emph{backward} if for $1\leq i+1<j \leq t$, $x_jx_i$ is an arc. Two or more paths are \emph{independent} if none of them contains an inner vertex of another. The \emph{length} of $P$ is the number of its arcs. A \emph{complete digraph} on $k$ vertices, denoted by $\mathop{K_k} \limits ^{\longleftrightarrow}$, is a digraph in which every pair of vertices is connected by an arc in each direction. As usual, a \emph{subdivision} of $\mathop{K_k} \limits ^{\longleftrightarrow}$ if a collection of $2\binom{k}{2}$ pairwise internally vertex disjoint paths joining every ordered pair of vertices in $\mathop{K_k} \limits ^{\longleftrightarrow}$. Further, for a positive integer $\ell$ we denote by $ T_{\ell}\mathop{K_k} \limits ^{\longleftrightarrow}$ a subdivision of $\mathop{K_k} \limits ^{\longleftrightarrow}$ where each arc is replaced by a path of length at most $\ell+1$. A subdivision is \emph{minimal} if all subdivision paths are backward.

We say that a digraph is \emph{regular} if all vertices have the same in-degree and out-degree. The unique acyclic tournament of order $n$ is the transitive tournament, which is denoted as $TT_n$. In this paper, we will often use the notation called \emph{blow-up}. Let $R$ be a digraph with vertex set $\{v_i: i\in [n]\}$, and let $G_1, \ldots,G_n$ be digraphs that are pairwise vertex disjoint. The blow-up of $R$ is the digraph $D$ with vertex set $V (G_1)\cup \ldots \cup V (G_n)$ and arc set $(\bigcup _{i=1}^n A(G_i))\cup \{g_ig_j: g_i\in V (G_i), g_j\in V (G_j)$, $v_iv_j\in A(R)\}$. 


\section{Counterexample}
Prior to a construction of the counterexample, it is necessary to introduce two key facts.
\begin{fact}\label{fac1}
\cite{Lichiardopol(2008)} Every regular tournament of order $n$ is $n/3$-connected.
\end{fact}

Let $l$ be a positive integer. Suppose $D$ is a digraph with a vertex $v$, and  a subset $S$. We use $(v, S, l,D)$ to represent that there are $l$ paths from $v$ to $S$ which only intersect at $v$ in $D$. Symmetrically, use $( S,v, l,D)$ to represent that there are $l$ paths from $S$ to $v$ which only intersect at $v$ in $D$.
\begin{fact}\label{fac2}
Let $D_1$ be a subdigraph of the digraph $D$, and let $S\subseteq V(D)$ be a subset of order at most $2k-1$ such that there is an $(x,y)$-path in $D\setminus S$, for any two vertices $x,y\in D_1\setminus S$. Given a vertex $v\in D\setminus D_1$, if there is $(v,V(D_1),2k,D)$ and $(V(D_1),v,2k,D)$ for $v$, then $D\setminus S$ contains a $(u,w)$-path, for any two vertices $u,w\in (D_1\cup v)\setminus S$.
\end{fact}
\begin{proof}
Fix two vertices $u,w\in (D_1\cup v)\setminus S$. By the assumption of $D_1$, we may suppose one of the vertices $u,w$ is $v$, say $u=v$. As a result of the assumption that there is still a path from $v$ to a vertex $x'\in D_1\setminus S$. Note that there is an $(x',w)$-path in $D\setminus S$. Hence, the $(u,w)$-path exists in $D\setminus S $. The proof for $w=v$ is analogous.
\end{proof}

The core of the proof of Theorem \(\ref{the1}\) lies in constructing counterexamples, which is accomplished by combining the subtournaments \(D_1\), \(D_2\) to obtain the final counterexample \(D\).

\begin{proof}\textbf{of Theorem \ref{the1}}
Let $m \geq 10k$ be an odd number. Firstly, we construct a tournament $D_1$ with  vertex set $V(D_1)=W\cup S\cup X'\cup Y$, where $|W|=m-2k+1$, $|S|=k-1$, $X'=$ $\left\{x'_1,  \ldots, x'_k\right\}$, and $ Y=\left\{y_1,  \ldots, y_k\right\}$. For convenience, we set $X'_2=X'\setminus \{x'_1\}$, and set $Y_2=Y\setminus \{y_1\}$. The arcs in $D_1$   satisfy the following properties (see Fig. \ref {fig1}):
\begin{itemize}
        \item $D_1\langle W\cup S\cup \{x'_1\}\cup Y_2\rangle$ forms a regular tournament, and $D_1\langle X'_2\rangle$ forms a tournament;
        \item $ W\cup \{x'_1\}\cup Y_2\rightarrow X'_2$, $X'_2 \rightarrow \{y_1\} \cup S$, and $ S\rightarrow y_1 \rightarrow  \{x'_1\}\cup W\cup Y_2$.
\end{itemize}

\begin{property}\label{prop1}
There is no $k$ disjoint $(x'_{\pi(i)},y_i)$-path in $D_1$, $i\in [k]$, where $\pi$ is any mapping of $[k]$ that satisfies $\pi(1)=1$.
\end{property}
\begin{proof}
As a result of our construction, any $(x'_{\pi(i)},y_i)$-path in $D_1$ that  internally avoids $X'\cup Y$ must use a vertex in $S$. Hence, there is no $k$ disjoint $(x'_{\pi(i)},y_i)$-paths, as $|S|=k-1$.
\end{proof}

Then we construct another tournament $D_2$ with $V(D_2)=\bigcup^k_{i=1}V(P^i) $. The arcs within $D_2$ meet the following properties (see Fig. \ref{fig2}):
\begin{itemize}
\item For each $i\in [k]$, the $P^i$ is a tournament with order $m$. Further, $P^i$ is a backward path, for all $i\in [k-1]$, and $P^k$ is regular. Namely, each $P^i$ has a Hamilton path. To simplify, we also use $P^i =p^i_1\cdots p^i_m$ to denote the Hamilton path of $P^i$, for all $i\in [k]$.
\item We construct a matching with $k-1$ arcs from $V(P^k)$ to $\bigcup^{k-1}_{i=1}V(P^i)$, i.e., set $E_1=\{p^k_{i} p^{i}_1\ |\ i\in [k-1]\}$. Additionally, we also construct $D_2$ so that $P^1$ cannot dominate any vertices in  $V(P^2[p^2_{m/2+1}, p^2_{m}]) $. Then we set $E_2=\{p^2_jp^1_l\ |\ j\in [m/2+1,m], l\in [m]\}$. All arcs except $E_1\cup E_2$ are directed such that $V(P^i)\rightarrow V(P^j)$, if $i<j$.
\end{itemize}

\begin{figure}[H]
\centering    
   \includegraphics[scale=0.4]{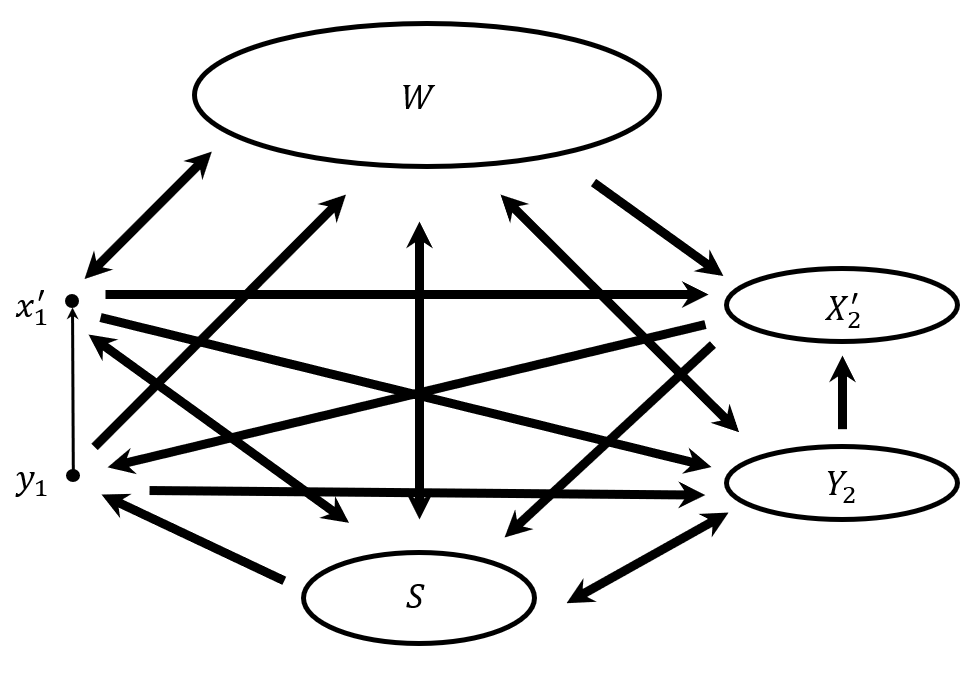}   
\caption{The structure of tournament $D_1$. The double arrow indicates that there may be arcs in opposite directions between the two sets, and the single arrow indicates that the direction of all arcs is the same as that of the arrow.}
\label{fig1}
\end{figure}

It can be easily observed that $N^-_{D_2}(P^1)=\{p^k_1\}\cup$ $ V(P^2[p^2_{m/2+1}$, $ p^2_{m}])  $ and $N^-_{D_2}(P^2[p^2_{m/2+1}$, $ p^2_{m}])$ $=\{p^2_{m/2}\}$ through our construction. Furthermore, we have the following property.
\begin{property}\label{proper2}
The digraph $D_2$ satisfies the following properties.
\begin{itemize}
\item[(i)] Any path from a vertex in $V(D_2\setminus P^2[p^2_{m/2+1},p^2_m])$ to a vertex $v$ on $P^2[p^2_{m/2+1},p^2_m]$ in $D_2$ must go through the vertex $p^2_{m/2}$ and then follow the path $P^2$ to $v$.
\item[(ii)] The path from $V(D_2\setminus P^1)$ to the vertex $p^1_m$, which avoids the vertices in $P^2[p^2_{m/2+1},p^2_m]$, must use the arc $p^k_1p^1_1$, and follow the path $P^1$ to reach the vertex $p^1_m$ in $D_2$.
\end{itemize}
\end{property}

Finally, let $D$ be a tournament on vertex set $V(D)=V(D_1) \cup V(D_2) \cup X$, where $D_1$ and $D_2$ are constructed as above and $X=\{x_1,x_2,\ldots, x_k\}$. Then $D$ has $(k+1)m+2k $ vertices. The arcs of $D$ are listed in the following way (see Fig. \ref{fig2}).
\begin{itemize}
  \item $D\langle X \rangle$ forms a tournament.
  \item $V(D_1)\setminus Y\rightarrow X$ except for the arc $ x_1x'_1$. Additionally, all arcs are directed from $X$ to $Y$ except for arcs between $x_i$ and $y_i$, for each $i \in [k]$.
  \item $N^+_{D_2}(x_i)=V(P^i)\setminus \{p^k_{k}, p^2_{m/2+1}, \ldots, p^2_{m}\}$, and $N^-_{D_2}(x_i)=V(D_2) \setminus N^+_{D_2}(x_i)$, for each $i\in [k]$, which ensures the large semidegree of each $x_i\in X$.
  \item We construct a matching with $k$ arcs from $D_2$ to $D_1$, i.e., set $E_3=\{p^i_mx'_i\ |\ i\in [k]\}$. Set $E_4=\{p^2_jy_1\ |\ j\in [m/2+1,m]\}$, which ensures that $y_1$ has sufficiently large in-degree. Then all arcs except $E_3\cup E_4$ are directed from $D_1$ to $D_2$.
\end{itemize}

\begin{figure}[H]
\centering    
   \includegraphics[scale=0.3]{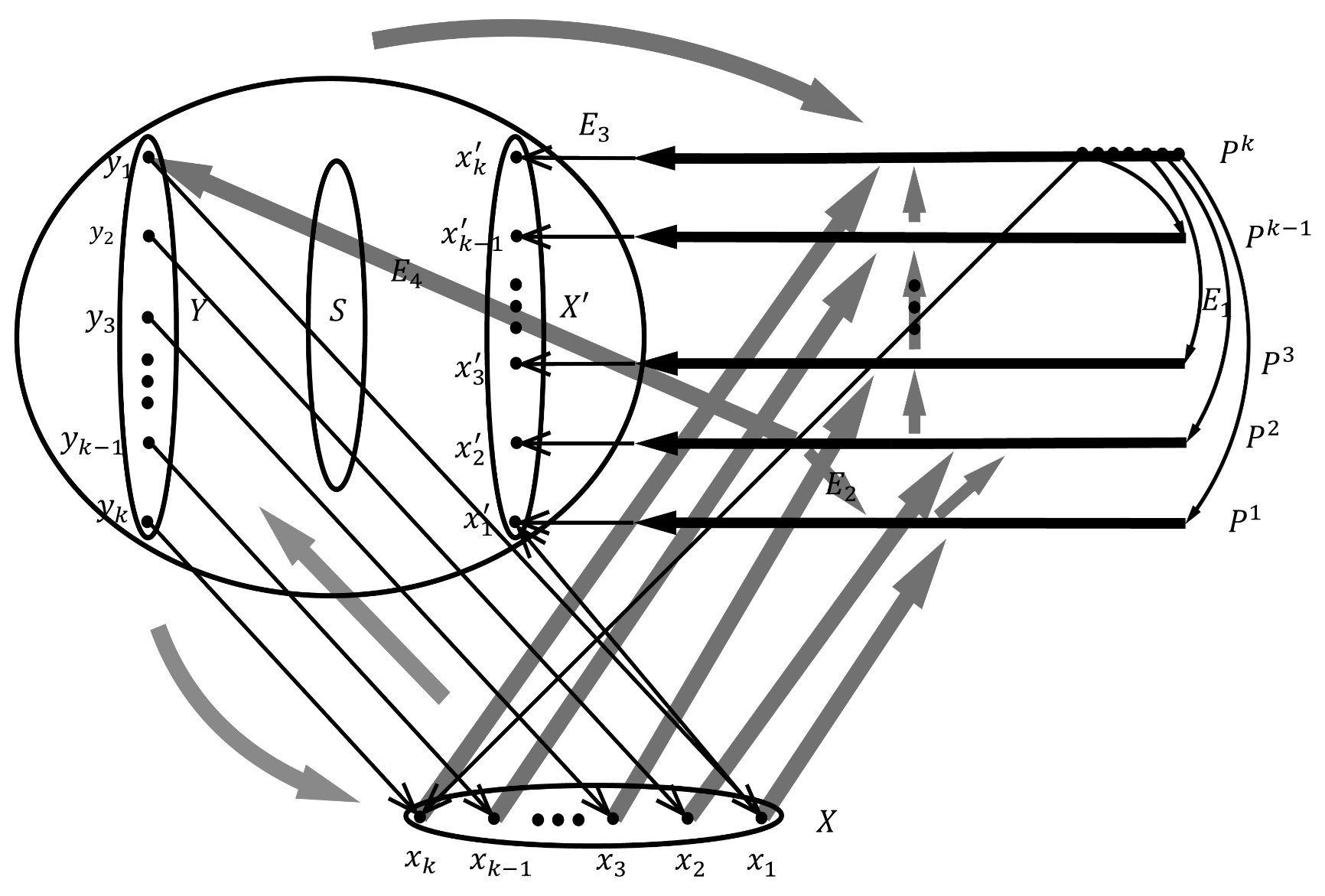}   
\caption{The counterexample $D$.}
\label{fig2}
\end{figure}

Claim \ref{cla1} validates that $D$ satisfies the conditions in Theorem \ref{the1}. Subsequently, Claim \ref{cla2} demonstrates that $D$ is not $k$-linked, thereby accomplishing the proof of the theorem.
\begin{claim}\label{cla1}
The tournament $D$ satisfies the following statements.
 \begin{itemize}
   \item[(P1)]  $\delta^0(D)\geq \lfloor m/2\rfloor= \lfloor \dfrac{|D|-2k}{2k+2}\rfloor$.
   \item[(P2)] $D$ is $2k$-connected.
 \end{itemize}
\end{claim}

\begin{proof}
It is straightforward to verify that $\delta^0(D)\geq \lfloor m/2\rfloor$. In order to confirm the connectivity of $D$, we need to show that for any subset $Z\subseteq V(D)$ of order at most $2k-1$ and any vertices $x,y \in D\setminus Z$, there is an $(x,y)$-path in $D\setminus Z$, fix such $Z$ in $D$ arbitrarily.

Firstly, we will show that there is an $(x,y)$-path in $D\setminus Z$, for any two vertices $x,y\in D_1\setminus Z$. If $x,y\in D[W\cup S\cup \{x'_1\}\cup Y_2]$, then the $(x,y)$-path in $D\setminus Z$ exists, as $D[W\cup S\cup \{x'_1\}\cup Y_2]$ is $2k$-connected by Fact \ref{fac1}. It is not difficult to see that $(W\cup S,y_1,2k,D )$ holds as $W\rightarrow X'_2\cup X\setminus \{x_1\}\rightarrow y_1$ and $S\rightarrow y_1$. Recall that $ y_1\rightarrow W$, then if $x,y\in D_1\setminus X'_2$, there exists an $(x,y)$-path in $D\setminus Z$ by Fact \ref{fac2}. Consider each vertex $x'_i\in X'_2$, we have the property $(x'_i,\{x_1,y_1\}\cup S\cup Y_2,2k,D)$. This is because $X'_2\rightarrow \{y_1\}\cup S\cup X$, $x_1x'_1\in A(D)$ and there are $k-1$ disjoint arcs from $X\setminus x_1$ to $Y_2$. Thus, combined with $W\rightarrow X'_2$, if $x,y\in D_1\setminus Z$, then there is an $(x,y)$-path in $D\setminus Z$ by Fact \ref{fac2}.

Next, we consider the connectivity of $D_1\cup D_2$. Recall that $P^k$ is a regular tournament, then $P^k$ is $2k$-connected by Fact \ref{fac1}. Further, consider each vertex $w\in V(D_2)$, say $w\in V(P^i)$. To apply Fact \ref{fac2}, we need to show that $w$ satisfies $(w,X'\cup Y,2k,D)$. Because $V(P^j)\rightarrow V(P^k[p^k_{k+1},p^k_m])$ for all $j<k$ and $P^k$ is $2k$-connected, we obtain the property $(w,\{ p^k_m, p^k_{1},\ldots, p^k_{k}\}$, $k+1,D\langle P^k\cup w\rangle)$. Using $E_1$ and the paths $P^i, i\in [k]$, we get $k$ internally disjoint paths from $w$ to $X'$, which are disjoint from $Y$. Together with the arcs in $\{p^k_{k}x_i\}\cup \{wx_j, j\in [k]\setminus i\}$, $w$ possesses the property $ (w,X\cup X',2k,D)$. Recall that there are $k$ disjoint arcs from $X$ to $Y$. Thus we get the property $(w,X'\cup Y,2k,D)$. Also, recall that every vertex $w\in D_2$ has at least $2k$ in-neighbours in $D_1$, there exists an $(x,y)$-path if $ x,y\in (D_1\cup D_2)\setminus Z$.

Similarly, we have  $W\rightarrow X$ and $x_i\rightarrow V(P^i)\setminus \{p^k_{m-k+1}, p^2_{m/2+1}, \ldots, p^2_{m}\}$, this also means that there is an $(x,y)$-path in $D\setminus Z$, for $x,y\in D\setminus Z$. The arbitrariness of $Z$ tells us that $D$ is $2k$-connected.
\end{proof}

\begin{claim}\label{cla2}
$D$ contains no $k$ disjoint paths linking $x_i$ to $y_i$, $i\in [k].$
\end{claim}

\begin{proof}
 To prove this, let us assume that such $k$ disjoint paths exist, written as  $\mathcal{Q}=\{Q_i, i\in [k]\ |\ Q_i$ is an $(x_i,y_i)$-path $\}$. Set $Q'_i$ be a subpath of $Q_i$ with the same initial vertex, whose terminal vertex $r_i$ is the first vertex in $D_1$ along the path $Q_i$, for all $i\in [k]$. Our construction shows that the second vertex, say $x^+_i$, of $Q_i$, must lie in $P^i$, for all $i\in [2,k]$. Recall that all arcs from $D_2$ to $D_1$ are contained in $E_1$. This also means that
\begin{equation}\label{5}
\text{the subpaths $Q'_i, i\in [2,k],$ must use $k-1$ arcs in $E_3$.}
\end{equation}

We show that $V(Q'_1)\cap V(P^2[p^2_{m/2+1},p^2_m])=\emptyset$. It can be seen from $N^+_{D\setminus (X\cup Y)}(x_1)$ $\subseteq \{x'_1\}\cup V(P^1)$ that the conclusion holds, if $x'_1$ is the second vertex of $Q_1$. Thus assume that the second vertex of $Q_1$ lies in the path $P^1$. If $V(Q_1[x_1,r_1])\cap V(P^2[p^2_{m/2+1},p^2_m] )\neq \emptyset$, then the vertex $p^2_{m/2}\in V(Q_1[x_1,r_1] )$ by Property \ref{proper2} (i). Together with Property \ref{proper2} (i)-(ii), any path $Q'_i$ ($i\in [2,k]$) does not use the arcs in $\{p^1_mx'_1, p^2_mx'_2\}$ to avoid the subpath $Q_1[x_1,r_1]$, which contradicts (\ref{5}). Thus $V(Q_1[x_1,r_1]\cap P^2[p^2_{m/2+1},p^2_m])=\emptyset$. Further, it is easy to see that there must have an arc in $E_3\setminus p^2_mx'_2$ intersecting the subpath $Q_1[x_1,r_1]$. This leads  to $V(Q_1[r_1,y_1])\cap V(P^2[p^2_{m/2+1},p^2_m] )= \emptyset$. Since otherwise, any path $Q'_i$ ($i\in [2,k]$) will not utilize the arc $p^2_mx'_2$, nor the arc that intersects with the path $Q_1[x_1,r_1]$, by Property \ref{proper2} (i)-(ii), which contradicts (\ref{5}). In summary, $V(Q_1)\cap V(P^2[p^2_{m/2+1}, p^2_m])=\emptyset$.

In fact, the collection $\mathcal{Q}$ must be joined by $k$ disjoint $(x_i,x'_{\pi(i)})$-paths from $X$ to $X'$ and $k$ disjoint paths from $X'$ to $Y$. Moreover, $\pi(1)=1$. To see this, we consider two cases. If $x^+_1=x'_1$, then $x'_1\in Q_1$. Together with $(\ref{5})$, $\{Q'_i, i\in [2,k]\}$ has to use the $k-1$ arcs in $E_3\setminus \{p^1_m x'_1\}$. Namely, $\mathcal{Q}$ are joined by $k$ disjoint $(x_i,x'_{\pi(i)})$-paths from $X$ to $X'$ and $k$ disjoint paths from $X'$ to $Y$. In this case, $\pi(1)=1 $. Now assume $x^+_1\in V(P^1)$, then $\mathcal{Q}$ has to use the $k$ arcs in $E_3$. Any path $Q'_i$ ($i\neq 1$) reaches $D_1$ either passes through the arc $p^1_mx'_1$ or through an arc in $E_3\setminus \{p^1_m x'_1\}$. In the former case, to avoid $Q'_1$, the path $Q'_i$ needs to use some vertices in $P^2[p^2_{m/2+1},p^2_m]$ to reach a vertex on $P^1$, then the arc $p^2_m x'_2$ must not be on any subpath $Q'_i, i\in [k],$ this contradicts the fact that $\mathcal{Q}$ utilizes the $k$ arcs in $E_3$. Then we conclude that $\mathcal{Q}\setminus Q_1$ must be joined by $k$ disjoint paths from $X$ to $X'\setminus \{x'_1\}$ and $k$ disjoint paths from $X'\setminus \{x'_1\}$ to $Y$ by (\ref{5}). Also, the fact that $\mathcal{Q}$ utilizes the $k$ arcs in $E_3$ concludes that $\pi(1)=1$.

Notice that the $k$ paths from $X$ to $X'$ already use at least one vertex from each arc in $E_3$, combined with $V(Q_1)\cap V(P^2[p^2_{m/2+1}, p^2_m])=\emptyset$, the $k$ disjoint paths from $X'$ to $Y$ must not use the vertices in $D_2$. Therefore these $k$ paths must be contained in $D_1$, but Property \ref{prop1} yields that there are no $k$ disjoint $(x'_{\pi(i)},y_i)$-paths in $D_1$ if $\pi(1)=1$, a contradiction.
\end{proof}
This completes our proof.
\end{proof}

\section{Proof of Theorem \ref{cor1}. }
In this section, let $D$ be a $(2k+1)$-connected semicomplete digraph with minimum out-degree at least $10^7k^{4}$. Suppose $X=\left\{x_1, \ldots, x_k\right\}$ and $Y=\left\{y_1, \ldots, y_k\right\}$ are two disjoint sets of $V(D)$. Our goal is to find the disjoint paths from $x_i$ to $y_i$, for each $i\in [k]$. We first give an outline of the proof of Theorem \ref{cor1}.
\subsection{Sketch of proof }
 First, we can show that the out-neighbourhood of $X$ has the following structure. There are three families $\mathcal{F} = \{F_{1},\ldots,F_{l}\}$, $\mathcal{S} = \{S_{1},\ldots,S_{l}\}$ and $\mathcal{V}=\{V_{l+1},\ldots, V_{k}\}$, for some $l\leq k$, which satisfy the following statements.
\begin{itemize}
\item $\{F_{1},\ldots,F_{l}\}$ are pairwise disjoint $T_2\mathop{K_{20k}} \limits ^{\longleftrightarrow}$. And the branch vertex set of $F_{i}$ is $S_{i}$, for all $i\in [l]$.
  \item $S_{i}\subseteq N_D^+(x_{i})$, for each $i\in [l]$, and $\{S_{1},\ldots,S_{l}\}$ forms a `almost' blow-up of a Hamiltonian path (see Fig. \ref{3a}).
  \item $V_j\subseteq  N_D^+(x_{j})$ is a set of order at least $16(k+1)$, for each $j\in [l+1,k]$, and $\{V_{l+1},\ldots, V_{k}\}$ forms a blow-up of transitive tournament according to Lemma \ref{lemma4} (see Fig. \ref{3b}).
\end{itemize}

By Menger's theorem, we can find a collection $\mathcal{Q}$ of $k+1$ disjoint paths from $O$ to $Y'$ which avoids the vertices in $X$, where $O$ is contained in $ S_{l}$ or $V_{k}$, and $Y'$ is the union of $Y$ as well as a `special' vertex in $V_{l+1}$. Our next goal is to link $X$ to $O$, that is, we need find a collection of $k$ short disjoint paths in $(D\setminus V(\mathcal{Q}))\cup O$. To achieve this, we need to prove that $\mathcal{Q}$ can be modified so that `many' of the subdivision paths of each $F_i$ do not intersect $\mathcal{Q}$ by Lemma \ref{lema1}, and $\mathcal{Q}$ also does not intersect each set $V_{j}$ in too many vertices by the transitivity of family $\mathcal{V}$.

\subsection{Preliminary lemmas}

In this subsection, we prepare some lemmas needed for our proof. The basic tool for finding disjoint paths is Menger's Theorem.

\begin{lemma}\textbf{\emph{(}Menger's Theorem\emph{)}}\label{lem}
\cite{Bang(1990)} Suppose $D$ is a $k$-connected digraph with two disjoint sets $X$ and $Y$. Then $D$ contains $k$ disjoint paths from $X$ to $Y$.
\end{lemma}

To show Lemma \ref{lemma4}, we need the following lemmas (Lemmas \ref{lemma2}-\ref{lemma1}). Lemma \ref{lemma2} provides an `almost' regular subset as the set of branch vertices for the subdivision. In order to demonstrate the inclusion relationship between $V_i$ and $W_i$ in Lemma \ref{lemma4}, we transform the problem into finding a matching in a bipartite graph, and Lemma \ref{lemma1} gives the necessary and sufficient condition for the existence of a matching in bipartite graphs.

\begin{lemma}\label{lemma2}
\cite{Snyder(2021)} There exists a positive integer $c$ such that any semicomplete digraph $D$ of order $n\geq s$ admits a subset $A$ of order $s$ with $ d^+(w)\geq \frac{9n}{40}$, $d^-(w)\geq \frac{9n}{40}$, and $d^-(w)\in [c-5s,c+5s]$, for every vertex $w \in A$.
\end{lemma}

\begin{lemma}\label{lemma1}
\cite{ Diestel(2016)} A bipartite graph with two parts $A$ and $B$ contains a matching covering $A$ if and only if $|N(S)|\geq |S|$ for all $S\subseteq A$.
\end{lemma}


Lemma \ref{lemma4} and Lemma \ref{lema1} are our key lemmas. Lemma \ref{lemma4} gives a blow-up of transitive tournament, which is used to characterize the structure of the out-neighbourhood of $X$. Here, we say that a subdivision \emph{sits} on $U$ if its branch vertex set is a subset of $U$.
\begin{lemma}\label{lemma4}
Let $\alpha$ be an integer. Suppose $D$ is a semicomplete digraph whose vertex set $V(D) = \bigcup_{i=1}^{\alpha}  U_i$, where $U_1,\ldots,U_{\alpha} $ are pairwise disjoint subsets with $|U_i|= 10^6k^3$ for all $ i\in [\alpha]$. If there is no $T_2\mathop{K_{20k}} \limits ^{\longleftrightarrow}$ sitting in any $U_i$, then there exists a blow-up $TT_{\alpha}[V_1,\ldots,V_{\alpha}]$ of a transitive tournament, where $V_i\subseteq U_i$ and $|V_i|\geq 16(k+1)$, for all $i\in [{\alpha}]$.
\end{lemma}
\begin{proof}
First, we need to find a blow-up of $TT_{5{\alpha}}$, say $BT_{5{\alpha}}[G_1,\ldots,G_{5{\alpha}}]$, such that all $G_i$ are non-empty and $|G_i|/5\leq |G_j|\leq 5|G_i|$ for all $i,j\in [5{\alpha}]$. To state this, we run the following process. Start with $G_1= D$ and $G_i=\emptyset$ for all $i\in [2,5{\alpha}]$. As long as there exists a empty set $G_i$, take $h$ to be a subscript satisfying $|G_{h}|=\max_{i\in [5{\alpha}]}  |G_i|$, and then use the `splitting operation' on $G_h$ to increase the number of non-empty sets. In the following, we will find a subset $B_h$ with order at least $20k$ such that $B_h\subseteq U_i$ for some $i\in [{\alpha}]$. With $B_h$ as the set of branch vertices, we construct a partial subdivision $M_h$ of $ T_2\mathop{K_{20k}} \limits ^{\longleftrightarrow}$ with the maximum number of subdivision paths. By assumption that there is no $ T_2\mathop{K_{20k}} \limits ^{\longleftrightarrow}$ sitting on $B_h$, there are two distinct vertices $u,v\in B_h$ such that
\begin{equation}\label{10}
 \text{$D\setminus M_h$ contains no $(u,v)$-path.}
\end{equation}
For convenience, denote $U'= N^+_{G_h}(u)\setminus V(M_h)$ and $V' = N^-_{G_h}(v)\setminus V(M_h)$. Thus $U'\cap V'=\emptyset$ and $ V'\rightarrow U'$ by (\ref{10}). Update $G_{j+1}:=G_{j} $ for $j\in \{5{\alpha}-1,\ldots,h+1\}$ successively, and then update $G_h:=V',\ G_{h+1}:=U'$, we get a new sequence $(G_1,\ldots, G_{5{\alpha}})$, this completes the splitting operation on $G_h$. Indeed, the existence of $B_h$ can be shown as follows. By Lemma \ref{lemma2}, there exists $c\in \mathbb{N}$ such that $G_h$ contains a subset $A_h$ of order $20k\alpha$, and for every vertex $ w \in A_h$
\begin{equation}\label{8}
  d^+_{G_h}(w)\geq \frac{9|G_h|}{40},\  d^-_{G_h}(w)\geq \frac{9|G_h|}{40},
\end{equation}
 and
\begin{equation}\label{11}
d^-_{G_h}(w)\in [c-100k\alpha,c+100k\alpha].
\end{equation}
We can partition $G_h$ into ${\alpha}$ parts $\{U_i\cap G_h, i\in [{\alpha}]\}$. By the pigeonhole principle, there is a subset $B_h$ of $A_h$ with order at least $20k$ such that $B_h\subseteq U_i$ for some $i\in [\alpha]$.


\begin{claim}\label{claimmm1}
In each splitting operation on $G_h$, the following statements are established.
\begin{itemize}
\item $|G_h|\leq |U'\cup V'|+800k^2+200k\alpha$.
\item After performing up to $5{\alpha}-1$ splitting operations, $|G_{h}|\geq 10^5k^2$. Moreover, after performing $5{\alpha}-1$ splitting operations, we get a new sequence $(G_1$, $\ldots,$ $G_{5{\alpha}})$ such that each $G_i$ is non-empty. That is, the programme terminates after $5{\alpha}-1$ splitting operations.
\item The new sequence $(G_1$, $\ldots,$ $G_{5{\alpha}})$ satisfies $|G_i|/5\leq |G_j|\leq 5|G_i|$ for all non-empty sets $G_i,G_j$.
\end{itemize}  
\end{claim}
\begin{proof}
Firstly, we show that $U'\cup V'$ covers at least $|G_h|-800k^2-200k\alpha$ vertices of $G_h$ in each splitting operation on $G_h$. By a simple calculate, $|V(M_h)|\leq 2\binom{20k}{2}+20k =400k^2$. Inequality (\ref{11}) yields that
\begin{equation*}
\begin{aligned}
200k\alpha &\geq |N^-_{G_h}(u)|-|N^-_{G_h}(v)|\\
&=|N^-_{G_h}(u)\setminus N^-_{G_h}(v)|+|N^-_{G_h}(u)\cap N^-_{G_h}(v)|-|N^-_{G_h}(v)\setminus N^-_{G_h}(u)|-|N^-_{G_h}(v)\cap N^-_{G_h}(u)|\\
&=|N^-_{G_h}(u)\setminus N^-_{G_h}(v)|-|N^-_{G_h}(v)\setminus N^-_{G_h}(u)|.
\end{aligned}
\end{equation*}
 Further, $|N^-_{G_h}(u)\setminus N^-_{G_h}(v)|\leq |N^-_{G_h}(v)\setminus N^-_{G_h}(u)| + 200k\alpha \leq 400k^2+200k\alpha$. This can be used to bound the lower bound of $U'\cup V'$ as following. 
\begin{equation}\label{1}
\begin{aligned}
|U'\cup V'|&\geq |(N^+_{G_h}(u)\cup N^-_{G_h}(v))\setminus M_h|\\
&\geq |N^+_{G_h}(u)\cup N^-_{G_h}(v)|-400k^2\\
&\geq  |G_h|-|N^-_{G_h}(u)\setminus N^-_{G_h}(v)|-400k^2\\
&>|G_h|-800k^2-200k\alpha .
\end{aligned}
\end{equation}

The inequality (\ref{1}) implies that at most $800k^2+200k\alpha$ vertices are lost in each splitting operation. Hence, $\sum_{i=1}^{5{\alpha}} |G_i|\geq |D|-(800k^2+200k\alpha) (5{\alpha}-1)$ after using up to $5{\alpha}-1$ splitting operations. The maximality of $G_h$ and $|D|\geq 10^6k^3$ mean that $$|G_{h}|\geq \frac{|D|-(800k^2+200k\alpha) (5{\alpha}-1)}{5{\alpha}}\geq 10^5k^2,$$ after $5{\alpha}-1$ splitting operations.
The inequality (\ref{8}) yields that both $U'$ and $V'$ have at least $9|G_h|/40 -400k^2$ vertices in each splitting operation on $G_h$. Namely, $|U'|,|V'|\geq |G_h|/5$ by $|G_h|\geq 10^5k^2$. Hence, every time we use a splitting operation, the new sequence $(G_1,\ldots, G_{5{\alpha}})$ has one less non-empty sets. So we get a new sequence $(G_1$, $\ldots,$ $G_{5{\alpha}})$ such that each $G_i$ is non-empty after performing $5{\alpha}-1$ splitting operations. 

Together with the maximum of $h$, we obtain that $|G_i|/5\leq |G_j|\leq 5|G_i|$, for all non-empty $G_i,G_j$. Since otherwise, without loss of generality, assume that there are two disjoint sets $G_i,G_j$ satisfying $|G_j|> 5|G_i|$. This shows that $G_j$ is a larger set than the one where the splitting operation was used to obtain $G_i$. The construction process suggests that a splitting operation on $ G_j$ would be performed before reaching $G_i$. Therefore $ G_i$ and $ G_j$ will not appear in the same sequence, a contradiction.
\end{proof}

 Claim \ref{claimmm1} follows that we can obtain a desired $BT_{5{\alpha}}[G_1,\ldots,G_{5{\alpha}}]$. To estimate the upper bound for all $G_i$, without loss of generality, set $|G_1|=\max_{i\in [5{\alpha}]} |G_i|$. Claim \ref{claimmm1} concludes that
\begin{equation*}
\begin{aligned}
{\alpha}|U_i|\geq \sum_{i\in [5{\alpha}]} |G_i|=|G_1|+\sum_{i\in [2,5{\alpha}]} |G_i|\geq |G_1|+\sum_{i\in [2,5{\alpha}]} |G_1|/5,
\end{aligned}
\end{equation*}
we have
\begin{equation}\label{3}
\begin{aligned}
|G_i|\leq |G_1|\leq \frac{5{\alpha} |U_i|}{5{\alpha}+4}<|U_i|, \text{ for each }i\in [5{\alpha}].
\end{aligned}
\end{equation}

We now construct a bipartite graph $H$ with two parts $H'=\{u_i,i\in [{\alpha}]\}$ and $H''=\{g_j, j\in [5{\alpha}]\}$, where $u_i$ corresponds to $U_i$ and $g_i$ corresponds to $G_i$. Let $u_ig_j$ be an edge of $H$ if and only if $U_i$ and $G_j$ intersect at least $16(k+1)$ vertices. Our purpose is to try to find a matching covering $H'$ by Lemma \ref{lemma1}, say $\{u_ig_{\pi(i)},i\in [{\alpha}]\}$. If this can be achieve, then there exists the desired $ TT_{\alpha}[V_1,\ldots,V_{\alpha}]$ by letting $V_i=G_{\pi(i)}\cap U_i$. Suppose there exists a subset $S$ of $H'$ such that $|N_H(S)|<|S|$. This means that
$$\bigcup_{u_i\in S} U_i\subseteq (D\setminus \bigcup_{i\in [5{\alpha}]}G_i)\cup \bigcup_{g_i\in N_H(S)}G_i\cup \bigcup_{g_i\notin N_H(S)}G_i,$$
that is,
$$ \sum_{u_i\in S} |U_i|< (800k^2+200k\alpha)(5{\alpha}-1)+\sum_{g_i\in N_H(S)} |G_i|+|S|(16(k+1)-1)(5{\alpha}-|N_H(S)|),$$
so $$|U_i|\leq (|S|-|N_H(S)|)|U_i|<(4000k^2{\alpha}+1000k{\alpha}^2)+90\alpha {(k+1)}|S|<10^6k^3.$$
This contradicts $|U_i|= 10^6k^3$. Hence, $|N_H(S)|\geq |S|$ for all $S\subseteq H'$. By Lemma \ref{lemma1}, $G$ has a matching covering $H'$.
\end{proof}

Lemma \ref{lema1} shows that we can modify the paths to release some vertices in a subdivision. To simplify the concept, given a subdivision $F$ whose branch vertex set is $S$, the subdivision path from $a$ to $b$ in $F$ is written as $P^F_{ab}$, for each ordered pair of vertices $a, b\in S$.

\begin{lemma}\label{lema1}
Suppose that the digraph $D$ contains a subdivision $F$ of the form $ T_2\mathop{K_{20k}} \limits ^{\longleftrightarrow}$ whose branch vertex set is $S$, and let $O=\{o_i,i\in [k+1]\}$ and $Y'=\{y_i, i\in [k+1]\}$ be two disjoint $(k+1)$-sets with $Y'\subseteq V(D\setminus F)$. Assuming $D$ has a collection $ \mathcal{Q}$ of $k+1$ disjoint paths from $O$ to $Y'$, then there exists a collection $\widehat{\mathcal{Q}}$ of $k+1$ disjoint paths from $O$ to $Y'$ such that there is a set $S'\subseteq S\setminus V(\widehat{\mathcal{Q}})$ of order $17k$ satisfying that $V(P^F_{uv})\cap V(\widehat{\mathcal{Q}})=\emptyset$, for any two vertices $u,v\in S'$. Moreover, if $ O\subseteq S$, the collection $\widehat{\mathcal{Q}}$ is the set of $k+1$ disjoint paths from some set $O'\subseteq S$ to $Y'$ such that $V(P^F_{uv})\cap V(\widehat{\mathcal{Q}})=\emptyset$, for any two vertices $u\in S\setminus O',v\in S$.
\end{lemma}
\begin{proof}
Our main idea is to use $ \mathcal{Q}$ to construct two sets of disjoint subpaths $\mathcal{Q}'$ and $\mathcal{Q}'' $, then connect $\mathcal{Q}'$ to $\mathcal{Q}'' $ by using disjoint subdivision paths of $F$. To ensure that the subdivision paths used are disjoint from $\mathcal{Q}''$, we sometimes have to update $\mathcal{Q}''$ to a new set $\mathcal{Q}'''$ of disjoint subpaths. Finally, we use disjoint subdivision paths to connect $\mathcal{Q}'$ to $\mathcal{Q}'''$, resulting in $k+1$ disjoint paths from $O$ to $Y'$ as desired, namely $\widehat{\mathcal{Q}}$. Assume that each path $Q_i\in \mathcal{Q}$ is an $(o_i,y_i)$-path.

To construct the desired collections of subpaths, we give the following definition. For a collection $\mathcal{Q}'=\{Q'_i, i\in [k+1]\ |\ Q'_i$ is a subpath of $Q_i\in \mathcal{Q}$, which starts at $o_i$ and ends at a vertex, say $f_i\}$, we call a subpath $Q'_i$ \emph{pre-good} with respect to   $\mathcal{Q}'$ if it satisfies one of the following statements:
\begin{itemize}
  \item[($Q'$1)] $Q'_i=Q_i$. In this case, $f_i=y_i$.
  \item[($Q'$2)] The terminal vertex $f_i$ of $Q'_i$ lies to a subdivision path $P^F_{a_ib_i}$ of $F$, where $b_i$ is not a pre-useful vertex of the other path $Q'_j\in \mathcal{Q}'$, and $P^F_{a_ib_i}[f_i,b_i]\cap V(\mathcal{Q}')=\{f_i\}$. In this case, we call the branch vertex $b_i$ \emph{pre-useful vertex} of $Q'_i$, call the path $P^F_{a_ib_i}[f_i,b_i]$ \emph{pre-useful path} of $Q'_i$. (Note that $P^F_{a_ib_i}[f_i,b_i]$ is the vertex $b_i$ if $b_i$ is a branch vertex.)
\end{itemize}


Evidently, there exists a collection $\mathcal{Q}'$ such that every path in $\mathcal{Q}'$ is pre-useful with respect to   $\mathcal{Q}'$, since $\mathcal{Q}$ is one. We choose such a collection $\mathcal{Q}'$ such that $\sum_{i=1}^{k+1}|Q'_i|$ is minimum (see Fig. 3(a)). From now on, let $\widetilde{S}_1(\mathcal{Q}')$ be a set of all non-pre-useful vertices of $\mathcal{Q}'$.

Symmetrically, for a collection $\mathcal{Q}''=\{Q''_i, i\in [k+1]\ |\ Q''_i$ is a subpath of $Q_i\in \mathcal{Q}$ that starts at a vertex, say $t_i$, and ends at $y_i \}$, we call a subpath $Q''_i$ \emph{post-good} with respect to $\mathcal{Q}''$ if it satisfies that
\begin{itemize}
  \item[($Q''$1)] $Q''_i=Q_i$. In this case, $t_i=o_i$.
  \item[($Q''$2)] The initial vertex $t_i$ of $Q''_i$ lies in a subdivision path $P^F_{c_id_i}$ of $F$, and $P^F_{c_id_i}[c_i,t_i]\cap V(\mathcal{Q}'')=\{t_i\}$, where $c_i\in \widetilde{S}_1(\mathcal{Q}')$ is non-post-useful for the other path in $\mathcal{Q}''$ and $d_i\in \widetilde{S}_1(\mathcal{Q}')$. In this case, we call the vertex $c_i$ \emph{post-useful vertex} of $Q''_i$, call the path $P^F_{c_id_i}[c_i,t_i]$ \emph{post-useful path} of $Q''_i$.
\end{itemize}
Clearly, we can choose a collection $\mathcal{Q}''$ such that every path in $\mathcal{Q}''$ is post-useful with respect to $\mathcal{Q}''$. Subject to this, we choose such a collection $\mathcal{Q}''$ such that $\sum_{i=1}^{k+1}|Q''_i|$ is minimum (see Fig. 3(b)). It is easy to check that all pre-useful vertices and all post-useful vertices are distinct. Also, each path $Q'_i\in \mathcal{Q}'$ (resp., $Q''_i\in \mathcal{Q}''$) has at most one pre-useful (resp., post-useful) vertex.

\begin{claim}\label{clai1}
 For any $i\in [k+1]$, the following statements hold.
 \begin{itemize}
   \item[(i)] $Q'_i$ is disjoint from any subdivision path $P^F_{uv}$, where $v$ is a non pre-useful vertex of $\mathcal{Q}'$.
   \item[(ii)] $Q''_i$ is disjoint from any subdivision path $P^F_{uv}$, where $u$ and $v$ are neither pre-useful for $\mathcal{Q}'$ nor post-useful for $\mathcal{Q}''$.
 \end{itemize}
\end{claim}
\begin{proof}
To prove $(i)$, assume that there exist some paths in $\mathcal{Q}'$ that intersect a such subdivision path $P^F_{uv}$. Choose $Q'_i$ and $w$ such that $w\in V(Q'_i)\cap V(P^F_{uv})$, and $P^F_{uv}[w,v]$ does not intersect with any pre-good path in $ \mathcal{Q}'\setminus Q'_i$. Thus each path in $\{ Q'_1,\ldots, Q'_{i-1}$, $Q'_{i}[o_i,w]$, $ Q'_{i+1}, \ldots, Q'_{k+1}\}$ is pre-good, which contradicts the minimality of $\sum_{i=1}^{k+1}|Q'_i|$. The proof of $(ii)$ is analogous.
\end{proof}

So far we have found two sets of subpaths $ \mathcal{Q}'$ and $  \mathcal{Q}''$. Our idea is that we connect the pre-useful paths and the post-useful paths by using $k + 1$ disjoint subdivision paths, so that we can connect $ \mathcal{Q}'$ and $  \mathcal{Q}''$ to get $k+1$ disjoint paths from $O$ to $Y'$. By the definition of pre-good path, the paths in $\{Q'_i[o_i,f_i]P^F_{a_ib_i}[f_i,b_i], i\in [k+1]\}$ are pairwise disjoint. Also, the definition of post-good path suggests that the paths in $\{P^F_{c_id_i}[c_i,t_i]Q''_i[t_i,y_i], i\in [k+1]\}$ are pairwise disjoint. Combined with Claim \ref{clai1}, $Q'_i$ is disjoint from all post-useful paths of $\mathcal{Q}'' $, because the terminal of any post-useful path is non-pre-useful for $\mathcal{Q}'$. Clearly, any pre-useful path is disjoint from post-useful paths, and $V(Q'_i)\cap V( Q''_j) = \emptyset$ for all $i\neq j$. But note that $Q''_i$ may intersect $Q'_i$, or the pre-useful path of $Q'_j$, or a subdivision path with a pre-useful initial vertex, for some $ i, j\in [k+1]$. Therefore, we need to update $\mathcal{Q}''$ again. 

For convenience, let $I=\{i\in I\ |\ $ $Q'_i\neq Q_i\}$, and let $U_i=Q'_i[o_i,f_i]P^F_{a_ib_i}[f_i,b_i]$ $\cup$ $\bigcup_{g_i \in \widetilde{S}_1(\mathcal{Q}')} P^F_{b_ig_i}$, for each $i\in I$. For simplicity, we use $U_i[a,b]$ to denote the path from $a$ to $b$ in $U_i$, for all vertices $a,b\in V(U_i)$. For a collection $\mathcal{Q}'''=\{Q'''_i, i\in [k+1]\ |\ Q'''_i$ is a subpath of $Q''_i\in \mathcal{Q}''$, which starts at a vertex, say $t'_i$, and ends at  $y_i \}$, we call a subpath $Q'''_i$ \emph{acceptable} with respect to   $\mathcal{Q}'''$ if it satisfies one of the following statements.
\begin{itemize}
  \item[($Q'''$1)] $Q'''_i=Q''_i$. 
  \item[($Q'''$2)] The initial vertex $t'_i$ of $Q'''_i$ is a vertex of $U_j$, where $j$ is not an acceptable index of the other path $Q'''_j\in \mathcal{Q}'''$. And $V(U_j[o_j,t'_i])\cap V(\mathcal{Q}''')=\{t'_i\}$. In this case, we call the index $j$ \emph{acceptable index} of $Q'''_i$, and call the path $U_j[o_j,t'_i]\setminus t'_i$ \emph{helpful path} of $Q'''_i$.
\end{itemize}

Similarly, we obtain a collection $\mathcal{Q}'''$, where all paths in $\mathcal{Q}'''$ are acceptable with $\mathcal{Q}'''$. We choose such a collection $\mathcal{Q}'''$ with minimum $\sum_{i=1}^{k+1}|Q'''_i|$ (see Fig. 3(c)).
\begin{figure}[htbp]
\centering
\begin{subfigure}[b]{0.325\textwidth} 
    \includegraphics[width=\textwidth]{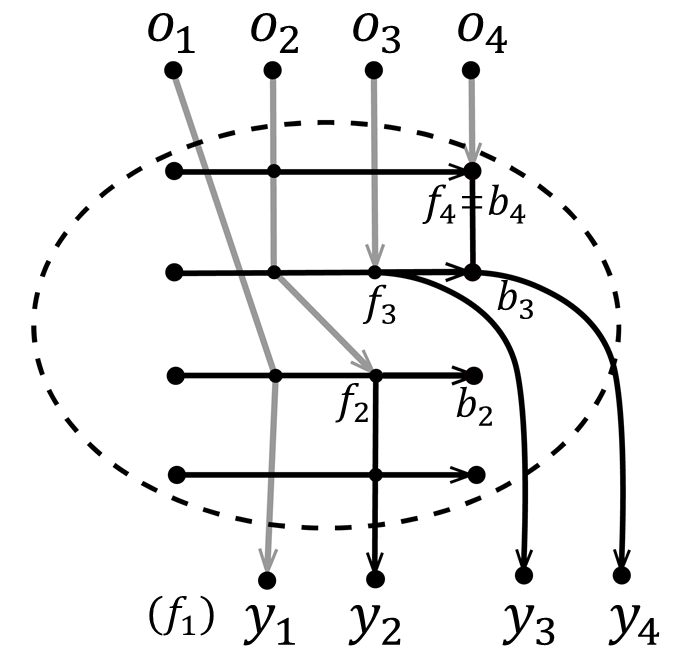}
    \caption{}
    \label{a}
\end{subfigure}
\begin{subfigure}[b]{0.32\textwidth}
    \includegraphics[width=\textwidth]{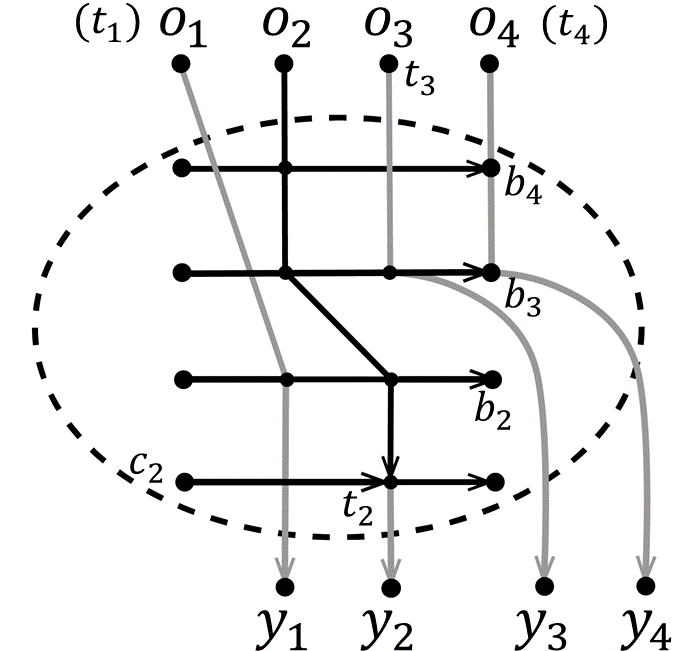}
    \caption{}
    \label{b}
\end{subfigure}
\begin{subfigure}[b]{0.275\textwidth}
    \includegraphics[width=\textwidth]{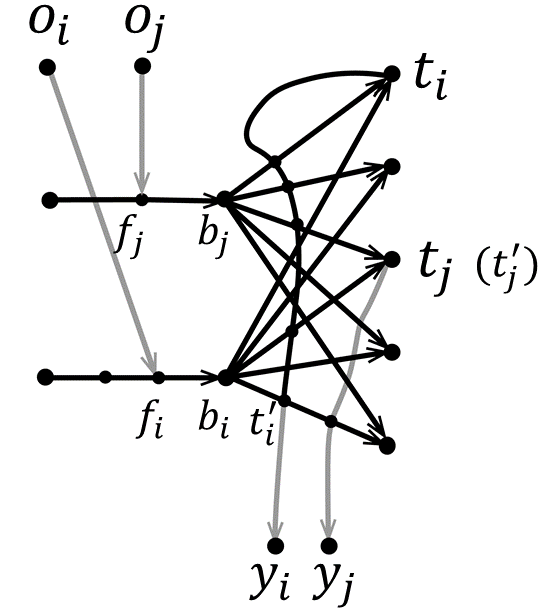}
    \caption{}
    \label{c}
\end{subfigure}
\caption{Figures (a), (b) and (c) represent the selection of $\mathcal{Q}'$, $\mathcal{Q}''$ and $\mathcal{Q}'''$ respectively. The black lines represent the subdivided paths and the initial paths $\mathcal{Q}$, while the gray lines represent the intercepted sub-paths. If the intersection vertices are not shown, it indicates that the two paths do not intersect.}
\end{figure}

One can directly obtain Claim \ref{claimm3} by the selection of $\mathcal{Q}'$ and $\mathcal{Q}''' $. 
\begin{claim}\label{claimm3}
The $\mathcal{Q}'''$ satisfies the following properties.
\begin{itemize}
\item[(i)] All helpful paths and all paths $Q'_i$, where $Q'_i\in \mathcal{Q}'$ satisfies $Q'_i=Q_i$, are pairwise disjoint.
\item[(ii)] $Q'''_i$ is disjoint from all helpful paths of $\mathcal{Q}'''$.
\item[(iii)] For any path $ Q'''_i\in \mathcal{Q}'''$ and for all indices $ j\in I$, $Q'''_i\cap U_j=\emptyset$ if $j$ is not the acceptable index of $\mathcal{Q}'''$.
\end{itemize}
\end{claim}

Now we connect the $k+1$ disjoint paths from $O$ to $Y'$ sequentially. Initially, set $K'=[k+1]$ to be the set of subscripts for paths in $\mathcal{Q}'$ that have not yet been used. And set $K''=[k+1]$ to be the set of subscripts for paths in $\mathcal{Q}'''$ that have not yet been used.
\begin{itemize}
  \item[(P1)] If there exists an index $i\in [k+1]$ such that $Q'_i=Q_i$, then let $\widehat{Q}_i=Q_i$, update $K'=K'\setminus i$, $K''=K''\setminus i$.
    \item[(P2)] If there exist two indexes $i\in K', j\in K''$ such that $i$ is the acceptable index of $Q'''_j$, let $\widehat{Q}_j=U_i[o_i,t'_j]Q'''_j[t'_j,y_j]$. Update $K'=K'\setminus i$, $K''=K''\setminus j$.
   \item[(P3)] For any two indexes $i\in K', j\in K''$, let $\widehat{Q}_j=Q'_i[o_i,f_i]P^F_{a_ib_i}[f_i,b_i]P^F_{b_ic_j}P^F_{c_jd_j}[c_j,t_j]Q''_j[t_j,y_j]$. Indeed, in this case, $Q'''_j=Q''_j$. Update $K'=K'\setminus i$, $K''=K''\setminus j$.
\end{itemize}

\begin{claim}
  $\widehat{\mathcal{Q}}=\{\widehat{Q}_i,i\in [k+1]\}$ is a collection of $k+1$ disjoint paths.
\end{claim}
\begin{proof}
Clearly, the paths selected by (P1) must be pairwise disjoint. It follows from Claim \ref{claimm3} (i)-(ii) that the paths constructed by (P2) must be pairwise disjoint. Arbitrarily select a path from (P1), say $\widehat{Q}_1=Q'_1=Q_1$, and a path from (P2), say $\widehat{Q}_j=U_i[o_i,t'_j]Q'''_j[t'_j,y_j]$, respectively. Indeed, $V(Q_1)\cap V( Q'_i[o_i,f_i]P^F_{a_ib_i}[f_i,b_i])=\emptyset$ for all $i\neq 1$, since $Q'_i$ is pre-good with respect to $\mathcal{Q}'$. By Claim \ref{clai1} (i), $Q_1$ is disjoint from the subdivision path with non-pre-useful terminal vertex. Thus, $V(Q_1)\cap V(U_i[o_i,t'_j])= \emptyset$. Apparently, $V(Q_1)\cap V(Q'''_j[t'_j,y_j])=\emptyset$ because $j\neq 1$. Thus $V(Q_1)\cap V(\widehat{Q}_j)=\emptyset$. Now, the paths constructed by (P1)-(P2) are pairwise disjoint.

We will assert that the paths constructed by (P3) must be pairwise disjoint. Suppose $\widehat{Q}_i=$  $Q'_{i'}[o_{i'},f_{i'}]P^F_{a_{i'}b_{i'}}[f_{i'},b_{i'}]P^F_{b_{i'}c_i}P^F_{c_id_i}[c_i,t_i]Q''_i[t_i,y_i]$ is a path generated by (P3), while $\widehat{Q}_{j}=Q'_{j'}[o_{j'}$, $f_{j'}]P^F_{a_{j'}b_{j'}}[f_{j'},b_{j'}]P^F_{b_{j'}c_j}P^F_{c_jd_j}[c_j,t_j]Q''_j[t_j,y_j]$ is another path produced by (P3). In (P3), none of the indexes in $K'$ is an acceptable index of $\mathcal{Q}'''$. This means that $V(Q''_j)\cap V(Q'_{i'}[o_{i'},f_{i'}]P^F_{a_{i'}b_{i'}}[f_{i'},b_{i'}]$ $P^F_{b_{i'}c_i})=\emptyset$. According to the definition of post-good, $Q''_j$ is disjoint from $P^F_{c_id_i}[c_i,t_i]$ $Q''_i[t_i,y_i].$ This yields that $ V(Q''_j)\cap V(\widehat{Q}_i)=\emptyset$. Analogously, $V(Q''_i)\cap V(\widehat{Q}_j)=\emptyset$. The definition of pre-good path and Claim \ref{clai1} (i) yield that $Q'_{j'}[o_{j'},f_{j'}]$ $P^F_{a_{j'}b_{j'}}[f_{j'},b_{j'}]P^F_{b_{j'}c_j}P^F_{c_jd_j}[c_j,t_j]$ is disjoint from $Q'_{i'}[o_{i'},f_{i'}]$ $P^F_{a_{i'}b_{i'}}[f_{i'},b_{i'}]P^F_{b_{i'}c_i}P^F_{c_id_i}[c_i,t_i]$. Consequently, $V(\widehat{Q}_i)\cap V(\widehat{Q}_j)=\emptyset$. The case where the path $\widehat{Q}_i $ is constructed by (P3) and $\widehat{Q}_{j}$ is a path produced by (P2) is analogous. Finally, we show that any path selected by (P1), say $\widehat{Q}_j=Q_j$, is disjoint from the path $\widehat{Q}_i$ which constructed by (P3). The definition of pre-good path and Claim \ref{clai1} (i) suggest that $V(Q_j)\cap V( Q'_{i'}[o_{i'},f_{i'}]P^F_{a_{i'}b_{i'}}[f_{i'},b_{i'}]P^F_{b_{i'}c_i}$ $P^F_{c_id_i}[c_i,t_i])=\emptyset$. Also, $V(Q_j)\cap V(Q''_i[t_i,y_i])=\emptyset$ as $i\in I$ (i.e. $i\neq j$). Thus $Q_j$ is disjoint from $\widehat{Q}_i $. In summary, the paths in $\widehat{\mathcal{Q}}$ are pairwise disjoint.
\end{proof}

By our construction, the collection $\widehat{\mathcal{Q}}$ occupies at most $2k+2$ vertices in $S$. Denote the set $S'$ to be the vertices in $S$ that are neither pre-useful for $\mathcal{Q}'$ nor post-useful for $\mathcal{Q}''$. It is easy to check that $S'$ has at least $20k-2(k+1)\geq 17k$ vertices.
\begin{claim}
  For any two vertices $u,v\in S'$, $P^F_{uv}\cap \widehat{ \mathcal{Q}}=\emptyset.$
\end{claim}
\begin{proof}
Here we give only those paths obtained by (P3) that do not intersect $P^F_{uv}$, since the rest of the cases are similar or even simpler.
For each path $\widehat{Q}_j$ obtained by (P3), which can be expressed as $\widehat{Q}_j=Q'_i[o_i,f_i]P^F_{a_ib_i}[f_i,b_i]P^F_{b_ic_j}P^F_{c_jd_j}[c_j,t_j]Q''_j[t_j,y_j]$, it should be noted that $V(Q'_i[o_i,f_i])\cap V(P^F_{uv})=\emptyset$ by Claim \ref{clai1} (i). Clearly, $P^F_{a_ib_i}[f_i,b_i]P^F_{b_ic_j}P^F_{c_jd_j}[c_j,t_j]$ is disjoint from $P^F_{uv}$. Claim \ref{clai1} (ii) implies that $Q''_j$ is disjoint from $P^F_{uv}$. Thereby, $V(P^F_{uv})\cap V(\widehat{ \mathcal{Q}})=\emptyset.$
\end{proof}

We have obtained the desired collection $\widehat{\mathcal{Q}}$ of $k+1$ disjoint paths from $O$ to $Y'$. Moreover, $S'$ is a subset of $S\setminus V(\widehat{\mathcal{Q}})$ with order $17k$, such that $P^F_{uv}$ is disjoint from $V(\widehat{\mathcal{Q}})$, for any two vertices of $u,v\in S'$. To prove ``moreover'' part of this lemma, we construct $\widehat{\mathcal{Q}}$ by just choosing a collection $\mathcal{Q}''$ such that every path in $\mathcal{Q}''$ is post-useful with respect to $\mathcal{Q}''$. Subject to this, we choose such a collection $\mathcal{Q}''$ such that $\sum_{i=1}^{k+1}|Q''_i|$ is minimum. Then we set $\widehat{\mathcal{Q}}=\mathcal{Q}''$. Similarly, $P^F_{uv}$ is disjoint from $V(\widehat{\mathcal{Q}})$, for any two vertices $u\in S\setminus O',v\in S$.
\end{proof}

\subsection{Proof of Theorem \ref{cor1}.}
\begin{figure}[htbp]
	\centering
	\begin{minipage}{0.43\linewidth}
		\centering
		\includegraphics[width=0.95\linewidth]{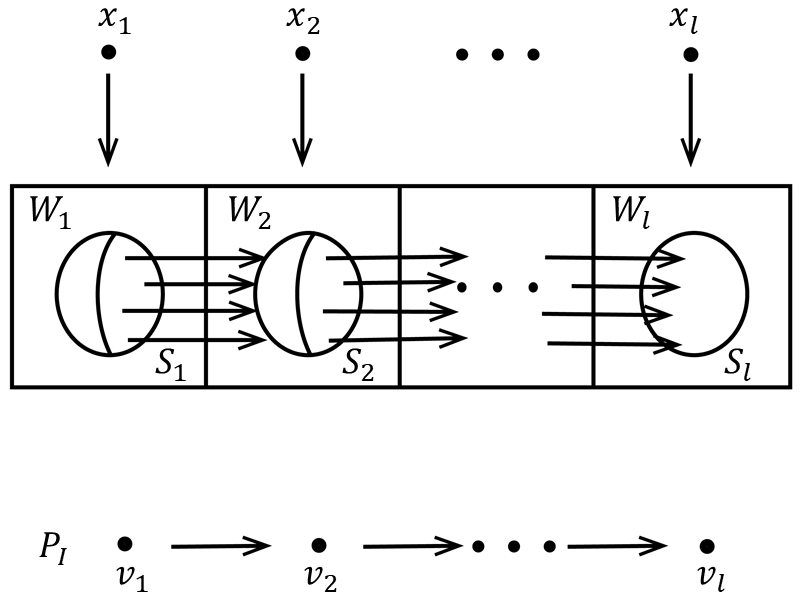}
		\caption{The construction of $H_I$.}
		\label{3a}
	\end{minipage}
	\centering
	\begin{minipage}{0.43\linewidth}
		\centering
		\includegraphics[width=0.95\linewidth]{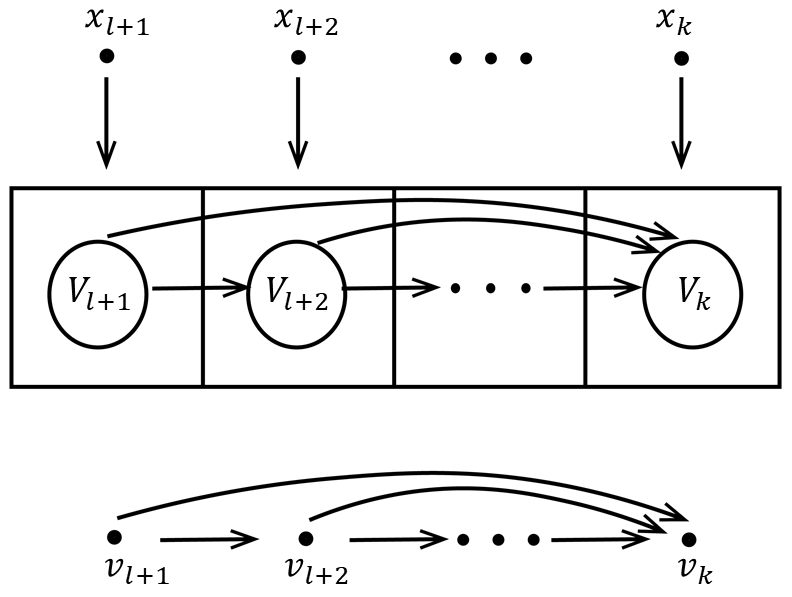}
		\caption{The construction of $\{V_{j}, j\in J\}$.}
		\label{3b}
	\end{minipage}
\end{figure}

Suppose that $D$, $X$ and $Y$ are given at the beginning of this section. First, we find the configuration in the out-neighbourhood of $X$. Because of the large minimum out-degree we can find a collection $\mathcal{W}$ of $k$ pairwise disjoint sets $W_i \subseteq N^{+}\left(x_i\right) \backslash (X\cup Y) $ with $\left|W_i\right|=(10^7 k^{4}-2k)/k $, $i \in [k]$. We select the index set $I$ as large as possible so that each $W_i$, $i\in I$, contains the branch vertices of $\mathop{T_2K_{20k}} \limits ^{\longleftrightarrow}$. Namely, we choose $I\subseteq [k]$ to be maximal such that
  \begin{itemize}
\item $\mathcal{F}$ contains $|I|$ mutually disjoint minimum subdivision $ \mathop{T_2K_{20k}} \limits ^{\longleftrightarrow}$;
\item for each $i\in I$, there is a subdivision $F_i\in \mathcal{F}$ sitting on the set $ W_i$.
\end{itemize}
Without loss of generality, assume $I= [l]$, and let $ J=[k]\setminus I$. For each $j\in J$, let $W'_j= W_j\setminus V(\mathcal{F}) $, then the order of $W'_j$ is at least
\begin{equation*}
\begin{aligned}
  \frac{10^7 k^{4}-2k}{k}-\left(20k+2\binom{20k}{2}\right)k & \geq \frac{10^7 k^{4}-2k}{k}-400k^3  \geq 10^6k^3.
\end{aligned}
\end{equation*}
 Apply Lemma \ref{lemma4} with $\alpha = |J|$ and $U_i\subseteq W'_i$ such that $|U_i|=10^6k^3$, $i\in J$, then there exist $|J|$ disjoint subsets $\{V_{j}, j\in J\}$ such that $V_{j'}\rightarrow V_{j''}$ for all $j'<j''$, which is denoted as $TT[V_{l+1},\ldots,V_k]$, and $V_{j}\subseteq W'_{j}$, $|V_{j}|\geq 16(k+1)$ for all $j\in J$ (see Fig. \ref{3b}).

Then define the auxiliary digraph $H_I$ with vertex set $\{v_i, i\in I\}$ as follows. For a pair $i_1, i_2 \in I$, $v_{i_1}v_{i_2}\in A(H_I)$ if and only if at least $|S_{i_1}|/4$ vertices in $S_{i_1}$ satisfy that each of them has at least $ \left|S_{i_2}\right| / 4$ out-neighbours in $S_{i_2}$. We claim that $H_I$ is a semicomplete digraph. Since $D$ is a semicomplete digraph, there are at least $|S_{i_1}||S_{i_2}|/2$ of the arcs between $S_{i_1}$ and $S_{i_2}$ that are in the same direction, and it may be assumed that these arcs are all from $S_{i_1}$ to $S_{i_2}$. If $S_{i_1}$ has less than $|S_{i_1}|/4$ vertices that has at least $ \left|S_{i_2}\right| / 4$ out-neighbours in $S_{i_2}$, then $S_{i_1}$ has at least $3|S_{i_1}|/4$ vertices with less than $ \left|S_{i_2}\right| / 4$ out-neighbours in $S_{i_2}$. Hence, the number of arcs from $S_{i_1}$ to $S_{i_2}$ is at most
$$ \dfrac{3|S_{i_1}|}{4}\cdot \dfrac{|S_{i_2}|}{4}+\dfrac{|S_{i_1}|}{4}\cdot |S_{i_2}|<\dfrac{|S_{i_1}||S_{i_2}|}{2},$$
a contradiction. So $H_I$ is a semicomplete digraph, and there is a Hamiltonian path in $H_I$. Without loss of generality, assume that it is this order that forms a Hamiltonian path, which is written as $P_I=v_1v_2v_3\cdots v_{l}$ (see Fig. \ref{3a}).  Note that one of $H_I$ and $J$ may be empty.

\textbf{Step 1. Select an appropriate set $O$ and a `special' vertex $y_{k+1}$, and then find $k + 1$ disjoint paths from $O$ to $Y\cup \{y_{k+1}\}$.}

Our aim now is to apply Menger's Theorem (Lemma \ref{lem}) to find $k+1$ disjoint paths from either $S_{l}$ or $V_{k}$ to $Y \cup\{y_{k+1}\}$ avoiding the vertices in $X$ if $J\neq \emptyset$, where $y_{k+1} \in V_{l+1}$. If $J=\emptyset$, then it is not necessary to take the vertex $y_{k+1}$; meanwhile, it is sufficient to take $k$ disjoint paths from $S_{l}$ to $Y$. In fact, we shall assume that $J\neq \emptyset$ since the proof of $J= \emptyset$ follows from the arguments for the case that $J\neq \emptyset$ (and is simpler). So we denote the set of the $k+1$ initial vertices of these paths by $O$. Define $V'_{k}=\{v\in V_{k}\ |\ d^+_{S_{l}} ( v)\geq \left|S_{l}\right|/ 2\}$, and define $V''_{k}=V_k\setminus V'_{k}$. The choice of $O$ and $y_{k+1}$ depends on the following statements:
\begin{itemize}
\item[(1)] If $P_I=\emptyset$ or $|V'_{k}|\geq |V_k|/2$, let $O \subset V'_{k}$ and let $y_{k+1} \in V_{l+1}$ be a vertex with maximum out-degree in $V_{l+1}$.
\item[(2)] If $J=\emptyset$ or $|V''_{k}|> |V_k|/2$, let $O \subset S_{l}$. Let $y_{k+1} \in V_{l+1}$ be a vertex with maximum out-degree in $V_{l+1}$ if $|J|\neq 1$. Let $y_{k+1} \in V''_{k}$ be a vertex with maximum out-neighbours in $V''_{k}$ if $|J|= 1$.
\end{itemize}

Since $D$ is $(2k+1)$-connected, there exists a collection $\mathcal{Q}$ of $k+1$ disjoint paths from $ O$ to $Y \cup\{y_{k+1}\}$ after removing the vertices in $X$, by Menger's Theorem. Choose $\mathcal{Q}$ with the minimum number of vertices. This also means that only the initial vertices of paths in $\mathcal{Q}$ are included in the corresponding set $S_{l}$ or $V_{k}$. Without loss of generality, assume that each $Q_{i}\in \mathcal{Q}$ is an $(o_{i},y_{i})$-path, where $i\in [k+1]$ and $o_i\in O$. We refer to the path in $\mathcal{Q}$ ending at a vertex that is not in $Y$ as the \emph{special path} in $\mathcal{Q}$, and we call the terminal vertex that is not in $Y$ the \emph{special vertex}. We denote that a vertex $s$ is $\mathcal{Q}$-free if $s\notin V(\mathcal{Q})$.

\textbf{Step 2. Adjust $\mathcal{Q}$ so that there are enough free vertices in $\{S_i,i\in I\}$ as well as in $\{V_j,j\in J\}$.}

 Using Lemma \ref{lema1} with each element $F_{i}\in \mathcal{F}$, we get a collection $\widehat{\mathcal{Q}}$ of $k+1$ disjoint paths from $O$ to $Y\cup \{y_{k+1}\}$ such that, for all $i\in I$, there exists a set $S'_{i}\subseteq S_{i}\setminus V(\widehat{\mathcal{Q}})$ of order $17k$ satisfying that
\begin{equation}\label{6}
\begin{aligned}
  &\text{$V(P^{F_{i}}_{uv})\cap V(\widehat{\mathcal{Q}})=\emptyset$, for any two vertices of $u,v\in S'_{i}$.}
\end{aligned}
\end{equation}
Specially, if $ O\subseteq S_{l}$, the collection $\widehat{\mathcal{Q}}$ is the set of $k+1$ disjoint paths from some set $O'\subseteq S_{l}$ to $Y\cup \{y_{k+1}\}$ such that $V(P^{F_{l}}_{uv})\cap V(\widehat{\mathcal{Q}})=\emptyset$, for any two vertices $u\in S_{l}\setminus O',v\in S_{l}$.




Next, we release the vertices in $V_{l+1},V_{l+2},\ldots,V_{k}$ in this order, by constantly modifying the collection $\widehat{\mathcal{Q}}$. We will show that there exists a collection $\mathcal{Q}^*$ of $k+1$ disjoint paths from $O$ to the union of $Y$ and a special vertex such that
\begin{itemize}
\item $|V_{j}\setminus V(\mathcal{Q}^*)|\geq 2$, for all $j\in [l+1,k-1]$,
\item $|V''_{k}\setminus V(\mathcal{Q}^*)|\geq 2$, if $|V''_{k}|\geq |V_k|/2$.
\end{itemize}
Note that there is nothing to prove if $t=k$ and $|V''_{k}|\leq |V_k|/2 $. Suppose that $|V''_{k}|> |V_k|/2 $. Without loss of generality, we assume that $\widehat{ \mathcal{Q}}$ satisfies $|V_{j}\setminus V(\mathcal{Q}^*)|\geq 2$, for all $j\in [l+1,t-1]$, for some $t< k$. Then we will construct a new collection $\widehat{ \mathcal{Q}}'$ to release enough vertices in $V_t$ as follows.



 Let $\widehat{ \mathcal{Q}}'=\widehat{ \mathcal{Q}}$, if
$$\begin{cases} |V_{t}\setminus V(\widehat{ \mathcal{Q}})|\geq 2,   \text{\ \  if $t\neq k$},\\ |V''_{k}\setminus V(\widehat{ \mathcal{Q}})|\geq 2, \text{\ \  if } t=k.
\end{cases}$$
Otherwise, we construct $\widehat{ \mathcal{Q}}'$ by means of path adjustments. First of all, we have the following inequalities.
 $$\begin{cases}
                |V(\widehat{ \mathcal{Q}})\cap V_{t}|\geq |V_{t}|- 1, & \mbox{if } t\neq k,\\
                |V(\widehat{ \mathcal{Q}})\cap V''_{k}|\geq |V''_{k}|-1, & \mbox{if } t=k.
              \end{cases}$$
 According to pigeonhole principle, there must have a path $Q \in \widehat{ \mathcal{Q}}$ satisfying 
 $$\begin{cases} |V(Q)\cap N_{V_{l+1}}^{+}(y_{k+1})|\geq \frac{16(k+1)-4}{4(k+1)}> 3, \text{ if }t=l+1,  \\|V(Q)\cap  V_t|\geq \frac{16(k+1)-1}{k+1}> 15,\text{ if }t\in [l+2,k-1],  \\|V(Q)\cap V''_{k}|\geq \frac{16(k+1)-2}{2(k+1)}> 7\text{, if }t=k.
\end{cases}$$
 Let $w_1, \ldots, w_{4}$ be the vertices in the intersection along $Q$. We obtain the $\widehat{ \mathcal{Q}}'$ by replacing $Q$ with $Q w_1$ if $Q$ is the special path in $\widehat{ \mathcal{Q}}$. And if $Q$ is not a special path in $\widehat{ \mathcal{Q}}$, assume $Q' \neq Q$ is the special path in $\widehat{ \mathcal{Q}}$ with terminal vertex $y_{k+1}'$. Then we will claim that $y_{k+1}'  w_{4}$ is an arc. When $t=l+1$, we obtain that $y_{k+1}'=y_{k+1}$ and $w_{4}\in N_{V_{l+1}}^{+}(y_{k+1}) $, so it is easy to check that $y_{k+1}  w_{4}$ is an arc in $D$. When $t\in [l+2,k]$, it is not difficult to find $y_{k+1}'\in V_j$ with $j<t$. Thus $y_{k+1}'  w_{4}$ is an arc in $D$ by the transitivity of $TT[V_{l+1},\ldots,V_k]$. Hence, we obtain desired collection $\widehat{ \mathcal{Q}}'$ by replacing $Q$ with $Q w_1$. And replace $Q'$ with the path following $Q'$ to $y_{k+1}' $, going along the arc $y_{k+1}' w_{4}$, then following $Q$ to its terminal in $Y$. In either case, the vertices $w_{2}, w_{3}$ are released, and $w_1$ is the special vertex of $\widehat{ \mathcal{Q}}'$. By induction, we can finally find the collection $\mathcal{Q}^*$ as desired. It is not hard to see that $ V(\mathcal{Q}^*)\subseteq V(\widehat{ \mathcal{Q}})$ by our construction. Without loss of generality, assume that each $Q^*_i\in \mathcal{Q}^*$ is an $(o_i,y_i)$-path, $i\in [k]$, this can be done by modifying the label of the vertices in $O$.




\textbf{Step 3. Find a collection $\mathcal{P}$ of $k$ disjoint paths linking $x_i$ to $o_i$ that do not use any vertices in $V(\mathcal{Q}^*)$.}

 It is easy to see that if we find such a collection, then we are done: for each $i\in [k]$ we can simply go from $x_i$ to $o_i$ via a path of $\mathcal{P}$ and then from $o_i$ to $y_i$ via a path of $\mathcal{Q}^*$, thus obtaining a path starting at $x_i$ and ending at $y_i$. The following claim is benefit for us to find disjoint paths in $ D\langle V(\mathcal{F})\rangle \setminus V(\mathcal{Q}^*)$. Recall that $S'_i$ of order $17k$ satisfies (\ref{6}).

\begin{claim}\label{claim7}
There are $2k$ disjoint paths from any $2k$-set $X'$ of $\bigcup_{i\leq l} S'_{i}$ to some $2k$-set of $S'_{l}$ in $D\setminus (X\cup V(\mathcal{Q}^*))$.
\end{claim}


\begin{proof}
It follows from the construction of the auxiliary digraph $H_I$ that at least $ \left|S_{i+1}\right| / 4$ vertices in $S_{i}$ satisfy that each of them has at least $ \left|S_{i+1}\right| / 4$ out-neighbours in $S_{i+1}$. Thus, consider the following two cases for each vertex $u\in X'$, say $u\in S'_{i}$. If $u$ has at least $5k$ out-neighbours in $S_{i+1} $, then $u$ has at least $5k-3k\geq 2k$ out-neighbours in $S'_{i+1} $. Let $ u^{+}\in S'_{i+1}$ be an out-neighbour of $u$. Otherwise, since $|S'_{i}|\geq |S_i|-3k$, there exist $\left|S_{i}\right| / 4-3k\geq 2k$ vertices in $S'_{i}$ such that each of them has at least $\left|S_{i+1}\right| / 4-3k\geq 2k$ out-neighbours in $S'_{i+1}$. Then $u$ can arrive to a vertex $u^{+}$ in $S'_{i}$ with at least $2k$ out-neighbours in $S'_{i+1}$ by a subdivision path in $F_{i}$, which avoids $V(\mathcal{Q}^*) $. Choose $ u^{++}\in S'_{i+1}$ to be an out-neighbour of $u^{+}$. Similar analysis for $u^{++}$, we can extend to obtain a minimal path from $u$ to a vertex in $S'_{l}$, say $u'$. Notice that a $(u,u')$-path use at most one vertex in $S'_{i_1}$ with $2k$ out-neighbours in $S'_{i_1+1}$, for each $i_1\in I$. Also, each $S'_{i_1}, i_1\in I,$ has at least $2k$ such vertices. So we can greedily find $2k$ disjoint paths from $X'$ to $S'_{l}$, which avoids $V(\mathcal{Q}^*) $.
\end{proof}

We then consider two cases depending on whether $O$ is in $V'_{k}$ or in $S_{l}$.

\textbf{Case 1. $O \subset V'_{k}$ }

Recall that for every $i \in J$ we have $V_{i} \subseteq N^{+}\left(x_{i}\right)$ and $V_{i} \rightarrow V_{i'}$ for each $i<i'$. For each $x_i\in X$, $i\in J$, choose a $\mathcal{Q}^*$-free vertex $x^+_i$ in $V_i$, so we can greedily link the vertices $x_i$ to $o_i$ by disjoint paths $x_ix^+_io_i$. On the other hand, to find a path from $x_i$ to $o_i$ for each $ i \in I$, we shall use the property that $o_i\in V'_{k}$ has $10k$ in-neighbours in $S_{l}$, so $o_i\in V'_{k}$ has $7k$ in-neighbours in $S'_{l}$. Thus, we may find a matching from a set $Z=\{z_i$, $ i \in I\}$ of $S'_{l}$ to $\{o_i, i \in I\}$ such that $z_i \rightarrow o_i$. 
For each $i\in I$, because $S'_{i}\setminus Z$ has at least $16k$ vertices, we can choose a matching from $\{x_i, i\in [k]\}$ to a subset $\{x^+_i, i\in [k]\}$ of $S'_{i}\setminus Z$ such that $x_i\rightarrow x^+_i$. Set $X^+_I=\{x^+_i, i\in I\}$. Apply Claim \ref{claim7}, there are disjoint paths from $X^+_I$ to some vertices in $S'_{l}\setminus Z$, say $Z'=\{z'_i$, $ i \in I\}$, in $D\setminus (X\cup V(\mathcal{Q}^*))$. Assume that such paths have the minimum number of vertices. Again using the subdivision paths of $F_{l}$, we can appropriately link $z'_i$ with $z_i$. Note that all vertices in $X^+_I\cup Z'\cup Z $ belong to the set $\bigcup_{i\in I}S'_i$, thus the paths are disjoint from $\mathcal{Q}^*$. Using the fact that $z_i\rightarrow o_i$, we obtain the desired paths from $x_i$ to $o_i$, where  $ i \in I$. Hence, we have found the required collection of disjoint paths linking $x_i$ to $y_i$, $i \in [k]$.

\textbf{Case 2. $O \subset S_{l}$}

\begin{figure}[htbp]
	\centering
	\begin{minipage}{0.8\linewidth}
		\centering
		\includegraphics[width=0.8\linewidth]{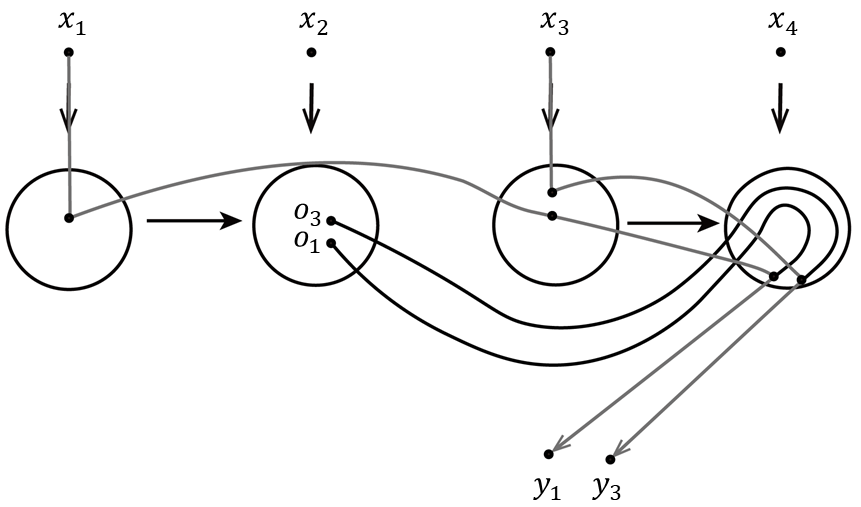}
		\caption{An example of Operation ($Q^* 1$)-($Q^* 2$).}
		\label{fig6}
	\end{minipage}
\end{figure}

 Recall that $O \cap S'_{l}=\emptyset$. In this case, $|V''_{k}|\geq |V_k|/2$. Set $J_{rem}$ to be a set of the indexes $i$ such that  $i\in J$ and we have not yet found the $(x_i,y_i)$-path. Initially, set $J_{rem}=J$, in the subsequent proof, the number of elements in $J_{rem}$ may decrease. From now on, we assert that $V''_{k}$ has $|J_{rem}|$ distinct $\mathcal{Q}^*$-free vertices by performing the following Operations ($Q^* 1$)-($Q^* 2$). On the contrary, assume that the number of $\mathcal{Q}^*$-free vertices in $V''_{k}$ is less than $|J_{rem}|$ now. Then we perform the following Operations ($Q^* 1$)-($Q^* 2$).
 \begin{itemize}
  \item[($Q^* 1$)] For all $i\in J_{rem}$, if there exists the $(o_i,y_i)$-path $Q^*_i$ in $\mathcal{Q}^*$ that intersects $V''_{k}$, let $v$ be the last intersecting vertex of $Q^*_i$ and $V''_{k}$. Then update $Q^*_i=x_ix^+_ivQ^*_i[v, y_i]$, where $x^+_i$ is a $\mathcal{Q}^*$-free vertices in $V_i$ and $x^+_i=v$ specially for $i=l$. Update $J_{rem}=J_{rem}\setminus \{i\}$ (see Fig. \ref{fig6}).  Otherwise, the operation terminates.
      \end{itemize}

 If $V''_{k}$ still has less than $|J_{rem}|$ distinct $\mathcal{Q}^*$-free vertices, run the ($Q^* 2$) operation. In the process of ($Q^* 2$), we may find some $(x_i,y_i)$-paths directly, $i\in I$. Here, denote $X_{used} $ to be a set of $\mathcal{Q}^*$-free vertices in the intersection of these $(x_i,y_i)$-paths and $V(\mathcal{F})$. $I_{rem}$ stores the indexes $i\in I$ and we have not yet found the $(x_i,y_i)$-path. Sometimes, we need to use $\mathcal{Q}^*$-free vertices in $V(\mathcal{F})$ to find $(x_i,y_i)$-path, where  $i\in J$, so we need to find a path from $x_i$ to $V(\mathcal{F})$ first. Then we let $X_{new}$ be a set of the terminal vertices of these paths. To ensure that we end up with enough free vertices in each $V_i$, we need to prevent $(Q^*2)$ from using the $\mathcal{Q}^*$-free vertices in $V''_l$ and using the $\mathcal{Q}^*$-free vertices in the same $V_i$ all the time. Therefore, we need to introduce a set $J_{unuse}$. Initially, set $X_{used}=\emptyset$, $I_{rem}=I$, $X_{new}=\emptyset$ and $J_{unuse}=\{k\}$.

\begin{itemize}
  \item[($Q^* 2$)] Choose $Q^*_i$ to be a path in $ \mathcal{Q}^*$ such that $|V(Q^*_i)\cap V''_{k}| $ is maximum for all $i\in I_{rem}$. Let $v$ be the last intersecting vertex of $Q^*_i$ and $V''_{k}$. Fix an index $i'\in  J_{rem}\setminus J_{unuse}$, we can choose $x^+_{i'}$ to be a $\mathcal{Q}^*$-free vertex in $V_{i'}$. This can be done since $|V_{i'}\setminus V(\mathcal{Q}^*)|\geq 2d^+_{H}(x_{i'})$. If $d^-_{S_{i}}(x^+_{i'})\geq 10k$, then $d^-_{S'_{i}}(x^+_{i'})\geq 7k$ as $|S'_i|\geq 17k$. Take $x^+_i$ to be a vertex in $S'_{i}\cap N^-(x^+_{i'}) $. Then update $Q^*_i=x_ix^+_ix^+_{i'}vQ^*_i[v, y_i]$. Since we have found the path from $x_i$ to $y_i$, update $I_{rem}=I_{rem}\setminus \{i\}$, $X_{used}=X_{used}\cup x^{+}_i$ and $J_{unuse}=J_{unuse}\cup \{i'\}$. Otherwise, take $x^{++}_{i'}$ to be an out-neighbour of $x^+_{i'}$ in $S'_{i}\setminus X_{new}$. Update $X_{new}=X_{new}\cup x^{++}_{i'}$, and $ J_{rem}=J_{rem}\setminus \{i'\}$ (see Fig. \ref{fig6}). Repeat this process until $J_{rem}\setminus J_{unuse}$ or $I_{rem}$ is empty.
\end{itemize}

\begin{claim}
$V''_{k}$ has at least $|J_{rem}|$ distinct $\mathcal{Q}^*$-free vertices.
\end{claim}
\begin{proof}
It is easy to see that Operation ($Q^* 1$) causes the value of $|J_{rem}|$ to decrease each time by 1. If $V''_{k}$ has no $|J_{rem}|$ distinct $\mathcal{Q}^*$-free vertices, then it means that $\bigcup_{i\in I_{rem}}V(Q^*_i)$ occupies at least $|V''_{k}|-k\geq 8(k+1)-k\geq k$ vertices in $V''_{k}$. Additionally, Operation ($Q^* 2$) indicates that either the value of $|J_{rem}|$ decreases by 1 or the number of $\mathcal{Q}^*$-free vertices increases by 2 without affecting the $\mathcal{Q}^*$-free vertices in $ J_{unuse}$ after using each Operation ($Q^* 2$). Operation ($Q^* 2$) can be executed until $J_{rem}\setminus J_{unuse}$ or $I_{rem}$ is empty. This indicates that the conclusion definitely holds if $J_{rem}\setminus J_{unuse}=\emptyset$. Note that the $(x_i,y_i)$-path obtained from ($Q^* 1$)-($Q^* 2$) intersects $V''_l$ at exactly one vertex. So if $I_{rem}$ is empty, together with the fact that Operation ($Q^* 1$) has already terminated, then $V''_l$ has at least $8(k+1)-k\geq k$ distinct $\mathcal{Q}^*$-free vertices as desired.
\end{proof}

For each $i\in J$, there exist at least two distinct $\mathcal{Q}^*$-free vertices in $ V_{i}$. There are at least $|J_{rem}|$ distinct $\mathcal{Q}^*$-free vertices in $ V''_{k}$, and every vertex in $ V''_{k}$ has at least $7k$ out-neighbours in $S'_{l}$. Thus, there exist $|J_{rem}|$ disjoint paths $P_i= x_ix^+_i x^{++}_iz_i$ for all $(i,j)\in J_{rem}$, where both $x^+_i\in V_{i}$ and $ x^{++}_i\in V''_{k}$ are $\mathcal{Q}^*$-free vertices, also $z_i$ is a vertex in $S'_{l}\cap N^+(x^{++}_i)$. Then update $X_{new}:=X_{new}\cup \{z_i, i\in J_{rem}\}$. For each $i\in I_{rem}$, choose $x^+_i$ to be an out-neighbour of $x_i$ in $S'_{i}\setminus X_{new}$. Then update $X_{new}:=X_{new}\cup \{x^+_i, i\in I_{rem}\}$. It is easy to check that $|X_{new}\cup X_{used}|\leq k$. By Claim \ref{claim7}, we can find disjoint paths from $X_{new}$ to a subset $Z'$ of $ S'_{l}$ in $D\setminus (X\cup V(\mathcal{Q}^*)\cup X_{used})$ with minimum number of vertices. Then use the subdivision paths in $F_{l}$ to link $Z'$ with corresponding vertices in $O$, which are disjoint from $V(\mathcal{Q}^*)$ by the ``moreover'' part of Lemma \ref{lema1}. Together with the disjoint paths that have been found, that is, the $(x_i,y_i)$-paths obtained directly by ($Q^* 1$)-($Q^* 2$), we obtain a collection of pairwise disjoint paths from $x_i$ to $y_i$ for every $i\in [k]$. This completes the proof of Theorem \ref{cor1}.

\section{Remark}

 Although we have determined that the connectivity of $2k + 1$ is tight and have also made efforts to reduce the bound of the minimum out-degree, we still haven't obtained a linear minimum out-degree, which conjectured by Gir\~{a}o, Snyder, and Popielarz in \cite{Girao(2021)}.
 \begin{conjecture}
\cite{Girao(2021)} There exists a constant $C >0$ such that every $(2k+1)$ connected tournament with minimum out-degree at least $Ck$ is $k$-linked.
 \end{conjecture}

Indeed, we can draw the following conclusion by considering the $H$-linkage problem, which is a generation of $k$-linkage. For two digraphs $H$ and $D$, an \emph{$H$-subdivision} $(f,g)$ in $D$ is a pair of bijections $f : V(H)\rightarrow V(D)$ and $g : A(H) \rightarrow \mathcal{P}(D)$ such that for every arc $uv\in A(H)$, $g(uv)$ is a path from $f(u)$ to $f (v)$, and distinct arcs map into internally vertex disjoint paths in $D$. Further, a digraph $D$ is \emph{$H$-linked} if every bijection $f : V(H) \rightarrow V(D)$ can be extended to an $H$-subdivision $(f, g)$ in $D$. Some results on $H$-linkage in graphs and digraphs also have been obtained previously. Please refer to \cite{Ferrara (2012), Ferrara (2006), Ferrara (2013), Gould(2006), kosto(2005), kostoc(2008),kostoc(20082),Wang(2024), Zhou(2024)}.  Theorem \ref{theorem} gives a connectivity for $H$-linked semicomplete digraphs. Here, we use $ H(h,m,\Delta,k_1,q)$ to denote a digraph $H$ of order $h$ with $m$ arcs and maximum out-degree $\Delta $ as well as $k_1$ vertices with in-degree zero, and let $q=\Delta\cdot \max \{h/2,m\}$. Notice that our aim is to find disjoint paths to replace the arcs in $H$, it is reasonable to assume that $H$ has no isolated vertices in this paper, i.e. $\delta (H)\geq 1$, where $\delta (H)=\min \{ d^+(v)+d^-(v)\ |\  v\in V(H)\}$.

\begin{theorem}\label{theorem}
Suppose $ H(h,m,\Delta,k_1,q)$ is a digraph with $\delta (H)\geq 1$. Then every $(k_1+m+1)$-connected semicomplete digraph $D$ with minimum out-degree at least $10^7q^{4}$ is $H$-linked.
\end{theorem}

 The proof of Theorem \ref{theorem} is similar to that of Theorem \ref{cor1}. Only the following improved version of Menger's theorem needs to be applied to find \(m + 1\) internally disjoint paths \(\mathcal{Q}\) after avoiding the vertices with in-degree 0 in \(H\). Meanwhile, it should be noted that we need to reserve $2d^+_H(x_i)$ distinct $\mathcal{Q}$-free vertices in each \(V_i\) during the proof. This is why we set the minimum out-degree to be \(10^7q^{4}\).
\begin{lemma}\textbf{\emph{(}Refining Menger's Theorem\emph{)}}
\cite{Bang(1990)} Let $D$ be a $k$-connected digraph. Let $x_1, \ldots, x_r$, $ y_1, \ldots, y_s$ be distinct vertices of $D$ and let $a_1,  \ldots, a_r$, $b_1,  \ldots, b_s$ be natural numbers such that
$$
\sum_{i=1}^r a_i=\sum_{j=1}^s b_j=k,
$$
Then $D$ contains $k$ internally disjoint paths $P_1,  \ldots, P_k$ with the property that precisely $a_i\left(b_j\right)$ of these start at $x_i$ (end at $\left.y_j\right)$.
\end{lemma}

 Thomassen \cite{Thomassen(19801)} proved that every 4-connected tournament has a Hamiltonian path from $x$ to $y$, for any two ordered vertices $x$ and $y$. It is therefore natural to ask the following.
 \begin{problem}
 Determine the optimal $f (k)$ such that every $f(k)$-connected tournament has $k$ disjoint paths $P_1,\ldots ,P_k$ such that $P_i$ is an $(x_i,y_i)$-path for all $i\in [k]$, and $\bigcup_{i\in [k]}V(P_i)=V(T)$, where $x_1,\ldots,x_k,y_1.\ldots,y_k$ are $2k$ distinct vertices.
 \end{problem}




\end{document}